\documentclass[openright,a4paper,american,amsart,11pt]{article}

\usepackage[latin1]{inputenc}

\usepackage{amsmath}

\usepackage{babel}

\usepackage{amsfonts}
\input xy
\xyoption{all}
\usepackage[dvips]{graphicx}
\usepackage{color}

\newtheorem{thm}{Theorem}[section]
\newtheorem{defi}[thm]{Definition}
\newtheorem{prop}[thm]{Proposition}
\newtheorem{lemma}[thm]{Lemma}
\newtheorem{cor}[thm]{Corollary}
\newcommand{\findemo}{\hfill
                     $\Box$ \vspace{1.5 ex}}
\newcommand{\RR}{\mathbb R}
\newcommand{\CC}{\mathbb C}
\newcommand{\ZZ}{\mathbb Z}

\begin{document}
\title{A VIRO THEOREM WITHOUT CONVEXITY HYPOTHESIS FOR TRIGONAL
  CURVES}
\author{Beno\^it Bertrand \and  Erwan Brugallé}
\date{}
\maketitle
\begin{abstract}
A cumbersome hypothesis for Viro patchworking of real algebraic curves
is the convexity of the given subdivision. It is an open question in
general to know whether the convexity is necessary. In the case of trigonal
curves we interpret Viro method in terms of \textit{dessins
d'enfants}. Gluing the dessins d'enfants in a coherent way we prove
that no convexity hypothesis is required to patchwork such
curves.
\end{abstract}

\setlength{\parindent}{0mm}

\section{Introduction}
In\footnote{Both authors are very grateful to the Max Planck
  Institute für Mathematik in Bonn for its financial support and
  excellent working conditions.} the late\footnote{\textit{2000
  Mathematics Subject 
  Classification : }14P25} seventies\footnote{Key words : topology of
  real algebraic curves, Viro method, dessin d'enfants, rational ruled
  surfaces} O. Ya. Viro
  invented a powerful method to 
construct real algebraic hypersurfaces with prescribed topology in
toric varieties. Since then it has become the main tool 
to prove the existence of certain topological types among families of
real algebraic varieties.

Classification of the real schemes realized by nonsingular curves of
degree $7$ in $\mathbb RP^2$ (\cite{V1}), smoothing of curves with
complicated singularities (\cite{V4}, \cite{Sh2}), maximal curvature
of plane algebraic curves (\cite{Lop}), construction of
counterexamples to the Ragsdale conjecture (\cite{I3}, \cite{I1},
\cite{Haa}, \cite{San} and \cite{Br2}) and construction of projective
hypersurfaces with big Betti numbers (\cite{B1}) are examples of
results using Viro method. It has also been used to prove maximality
results. A real algebraic variety $X$ is said to be maximal ($M$-variety)
if the sum of the Betti numbers with $\ZZ_2$ coefficients of its real
part is equal to the corresponding sum for its complexification
$(b_*(\RR X;\ZZ_2)=b_*(\CC X;\ZZ_2))$. Viro method has been very
useful to prove for instance the existence of projective
$M$-hypersurfaces of any degree in any dimension (\cite{IV1}) or
existence of asymptotically maximal families of hypersurfaces in any
projective toric variety (\cite{Ber1}).

Roughly speaking , the Viro method (for curves) can be described as
follows~: start with some real algebraic curves $(C_i)$ whose Newton
polygons are the $2$-cells of a polygonal subdivision $\tau$ of a convex
polygon $\Delta$.  Then, under some assumptions on the $C_i$'s and on
$\tau$, the Viro Theorem asserts that there exists a real algebraic
curve with Newton polygon $\Delta$ 
obtained as a topological gluing of the curves $C_i$.
One of
the hypothesis of Viro's theorem is that the subdivision $\tau$ should
be convex, i.e. there should exist a convex piecewise-linear affine
function $\lambda:\Delta\to\mathbb R$ whose domains of linearity are
exactly the Newton polygons of the curves $C_i$.
Notice that non-convex subdivisions do exist (see \cite{CoHe}),  an example is
depicted in Figure \ref{non convex subd}.

It is an open problem to know whether the convexity hypothesis is
necessary or not. More precisely it is not known whether a curve
patchworked from a non-convex subdivision is isotopic to a real
algebraic curve with Newton polygon $\Delta$ (see \cite{Loe} and
\cite{IS}). I. Itenberg and E. Shustin proved in \cite{IS} that in the
case of $\mathbb CP^2$ or of the geometrically ruled rational surfaces
$\Sigma_n$ there always exists a real pseudoholomorphic curve
isotopic to the patchworked one. 
Note that it is still unknown if there exists a real
nonsingular pseudoholomorphic curve in $\mathbb RP^2$
which is not isotopic to a real nonsingular algebraic curve of the
same degree. 

If the convexity hypothesis in the Viro method turned out
to be necessary in general it could lead to such examples of
``non-algebraic pseudoholomorphic curves''. 
On the other hand if the
convexity hypothesis could be removed it would simplify greatly the Viro
construction since it is often rather technical to check whether a given
subdivision is convex. Moreover this approach would give a deeper
insight or at least a new interpretation of the Viro method. 
This work is a first step in the study of the need of the convexity
hypothesis.  
\newline

Trigonal curves are those curves in $\Sigma_n$ whose Newton polygon is 
the triangle $\Delta_n$
with vertices $(0,0)$, $(0,3)$ and $(3n,0)$.
In this article we prove that we can get rid of the convexity hypothesis
for trigonal curves. Namely, any curve patchworked from $\Delta_n$ is
an algebraic trigonal curve in $\Sigma_n$ (with Newton polygon
$\Delta_n$). In fact we construct curves with prescribed positions
with respect to the natural pencil of lines of $\Sigma_n$. In the
convex case, the pencil of lines has been studied by  L. Lopez de
Medrano  in 
\cite{Lop}.
  
We want to construct a trigonal curve $C$ out of a trigonal patchwork.
The strategy is the following.  We first cut $C$ in several pieces
which are the intersection of $C$ with some portions of the pencil of
lines in $\Sigma_n$.  All these pieces are given by the subdivision of
$\Delta_n$ and will be algebraic. In \cite{IS}, Itenberg and Shustin
glue these pieces topologically in order to obtain a real
pseudoholomorphic curve. 

Since we deal with trigonal curve we are able to use a particular case
of "dessins d'enfants" which are called real rational graphs. This
technic was introduced in real algebraic geometry independently by S.
Yu. Orevkov in \cite{O3} and by S. Natanzon, B. Shapiro and A.
Vainshtein in \cite{Shap1}.  The position of a real trigonal curve
with respect to the pencil in $\Sigma_n$ is encoded by three linearly
dependent real polynomials (see \cite{O3} or section \ref{poly deg
  3}).  Real rational graphs give a necessary and sufficient condition
for the existence of three linearly dependent real polynomials in one
variable whose roots realize a given real arrangement. Via Riemann
existence theorem they provide an existence criterion for a given
trigonal curve.

In fact, we will use a slight generalization
of these objects. To glue algebraically the pieces of our curve we
indeed need to perform some surgery on dessins d'enfants enhanced with
the data of the sign of some characteristic polynomials. These are what
we call signed real rational graphs. Once glued together in a coherent way
the signed real rational graphs obtained for each piece yield a 
dessin d'enfant which corresponds to a curve with the required topology.

\vspace{2ex}
\textbf{Organization of this article}

In section \ref{defi}, we recall some facts about rational
geometrically ruled surfaces. In section \ref{state viro}, we explain
the patchwork construction and  give a patchwork theorem for
 real algebraic trigonal curves without any convexity assumption (Theorem
 \ref{viro}). We explain in section \ref{poly deg 3} how to encode the
 topology of a trigonal curve in a sign array. In section \ref{pencil}
 we 
define
 an order on the polygons of a patchwork
which is used in section \ref{sign patchwork} to extract a sign
array from a trigonal patchwork.
We state there our main result (Theorem \ref{non convex viro})
. The rest of the paper is devoted
to the proof of this theorem. The main tool of this proof, signed real
rational graphs, is defined in
section \ref{graph}. In section \ref{constr graph} we associate a
signed real
rational graph to each piece of the pencil of line given by the order
on the polygons of a patchwork. We glue all these graph in section
\ref{gluing graphs} and check in section \ref{proof} that the obtained
signed real rational graph corresponds to a real algebraic trigonal
curve which has the topology prescribed by the patchwork.

\vspace{2ex}
\textbf{Notation}

All  polynomials in two variables $C(X,Y)$ are
considered as polynomials in the variable $Y$ whose coefficients are
polynomials in the variable $X$. 
By the discriminant of $C(X,Y)$ 
we mean the  discriminant of $C(X,Y)$ with respect to the variable $Y$.

\vspace{2ex}
Suppose that a coordinate system of $\mathbb RP^1$ is fixed. The point
$[0:1]$ (resp. $[1:0]$) is denoted by $0$ (resp.$\infty$), and
the points $[x:y]$ such that $xy>0$ (resp. $xy<0$) are called
positive (resp. negative). The embedding of $\mathbb R$ into $\mathbb
RP^1$ given by $x\mapsto [x:1]$ induces an orientation of $\mathbb
RP^1$. Given two points $A$ and $B$ in $\mathbb RP^1$, the segment
$[A;B]$ is the connected component of $\mathbb RP^1\setminus
\{A,B\}$ oriented from $A$ to $B$ for this orientation. 

\vspace{2ex}
When there is no ambiguity, the two words ``curve'' and ``polynomial''
are used to designate either a polynomial $C(X,Y)$ or the curve
defined by this polynomial.

\section{Rational geometrically ruled surfaces}\label{defi}
\textit{The $n^{th}$ rational geometrically ruled surface}, denoted by
$\Sigma_n$, is the surface obtained by taking four copies of $\mathbb
C^2$ with coordinates $(x,y)$, $(x_2,y_2)$, $(x_3,y_3)$ 
and $(x_4,y_4)$, and by gluing them along $(\mathbb C^*)^2$ with the
identifications $(x_2,y_2)=(1/x,y/x^n)$, $(x_3,y_3)=(x,1/y)$ 
and $(x_4,y_4)=(1/x,x^n/y)$. Let us denote by $E$  (resp. $B$ and $F$)
the algebraic curve in $\Sigma_n$ defined by the equation $\{y_3=0\}$
(resp. $\{y=0\}$ and $\{x=0\}$). The coordinate system $(x,y)$ 
is called \textit{standard}.
The projection $\pi$~: $(x,y)\mapsto x$ on $\mathbb C^2$ defines a
$\mathbb C^1$-bundle structure on $\mathbb C^2$ which is extendable up to a 
$\mathbb CP^1$-bundle structure on $\Sigma_n$. 
The intersection numbers of $B$ and $F$ are respectively
$B\circ B=n$, $F\circ F=0$ and $B\circ F=1$. The surface
$\Sigma_n$ has a natural real structure induced by the complex
conjugation in $\mathbb C^2$, and 
the real part $\RR\Sigma_n$ of $\Sigma_n$ is a torus if $n$ is even and a Klein
bottle if $n$ is odd. The restriction of $\pi$ on $\mathbb R\Sigma_n$
defines a pencil of lines denoted by $\mathcal L$. 

 The group $H_2(\Sigma_n,\mathbb Z)$ is isomorphic to $\mathbb
 Z\times\mathbb Z$ and is generated by the classes of $B$ and
 $F$. Moreover, one has $E=B-nF$. An algebraic curve on $\Sigma_n$ is 
said to be of \textit{bidegree} $(k,l)$ if it realizes the homology
 class $kB+lF$ in $H_2(\Sigma_n,\mathbb Z)$. Its equation in
 $\Sigma_n\setminus E$ is
\begin{displaymath}
\sum_{i=0}^k a_{k-i}(X,Z)Y^i
\end{displaymath}
where $a_j(X,Z)$ is a homogeneous polynomial of degree $nj+l$.
A curve of bidegree
 $(3,0)$ is called \textit{a trigonal curve} on $\Sigma_n$. 

In the rational geometrically ruled surfaces, we study real curves up
to \textit{isotopy with respect to $\mathcal L$}. Two curves are said 
to be isotopic with respect to the fibration $\mathcal L$ if there
exists an isotopy of $\RR\Sigma_n$ which transforms the first curve to
the second one, and 
which maps each line of $\mathcal L$ 
to another line of $\mathcal L$.
In this paper, curves in a rational geometrically ruled surface are
depicted up to isotopy with respect to $\mathcal L$. 

\begin{defi}
An algebraic curve $C$ of bidegree $(k,l)$ in $\Sigma_n$ is said to be $\mathcal
L$-nonsingular if $C$ is nonsingular and if for any fiber $F$ in
$\Sigma_n$, the set $F\cap C$ contains at least $k-1$ points.
\end{defi}

\section{Patchworking real algebraic trigonal curves}\label{state viro}

Here, we explain the patchworking method for curves. The general patchworking
Theorem, which uses the convexity of the subdivision, can be found in
\cite{V1}, \cite{V4},   \cite{V2}, \cite{Ris} and \cite{IS2}. Our patchwork
Theorem, Theorem \ref{viro}, does not require this convexity
assumption, yet it works only for trigonal curves.

\subsection{Chart of a real polynomial}
This is the key notion of the Viro method. 
A chart of a real
polynomial $C(X,Y)$ is a way to draw the curve
defined by $C$ in $(\mathbb R^*)^2$ in the union of 4 symmetric copies of its Newton polygon. The Newton polygon of
 a polynomial 
$F$ will be denoted by $\Delta(F)$.

Given a  convex polygon $\Delta$ in $\mathbb R^2$ with vertices in
$\mathbb Z^2$,
one can define the 
 so called
 moment map from $(\mathbb R^*_+)^2$ to
the interior of $\Delta$ as follows
\begin{displaymath}
\begin{array}{cccc}
\mu_{\Delta}~:&(\mathbb R^*_+)^2&\to&I(\Delta)\\
   &   (x,y)&\mapsto& \frac{\sum_{(i,j)\in
   \Delta\cap\mathbb Z^2}x^{i}
  y^{j}.(i,j)}{\sum_{(i,j)\in \Delta\cap\mathbb Z^2}x^{i}y^{j}}
\end{array}.
\end{displaymath}
If $\dim(\Delta)=2$, then $\mu_\Delta$ is a diffeomorphism.
For a pair $(\epsilon_1,\epsilon_2)\in\{+,-\}^2$, we denote by
$s_{\epsilon_1,\epsilon_2}$ the symmetry of $\mathbb R^2$ given by
 $s_{\epsilon_1,\epsilon_2}(x,y)=(\epsilon_1x,\epsilon_2y)$. For a
convex polygon $\Delta$ in $(\mathbb R_+)^2$, we define $\Delta^*$ as
$\bigcup_{(\epsilon_1,\epsilon_2)\in \{+,-\}^2}s_{\epsilon_1,\epsilon_2}(\Delta)$.

\begin{defi}
Let $C(X,Y)$ be a real polynomial. The chart of $C(X,Y)$, denoted by
$Ch(C)$, is the closure of the set
\begin{displaymath}
\bigcup_{(\epsilon_1,\epsilon_2)\in \{0,1\}^2}
s_{\epsilon_1,\epsilon_2}\circ \mu_{\Delta(C)}\circ
s_{\epsilon_1,\epsilon_2}(\{(x,y)\in\mathbb R^*_{\epsilon_1}\times \mathbb
R^*_{\epsilon_2}| C(x,y)=0\})  
\end{displaymath}
\end{defi}
It is clear that $Ch(C)\subset\Delta^*$.

\subsection{Gluing of charts}

Let $C(X,Y)=\sum a_{i,j}X^iY^j$ be a polynomial and $\gamma$ be a face
of $\Delta(C)$.
\begin{defi}
The truncation of $C(X,Y)$ on $\gamma$ is the polynomial
$\sum_{(i,j)\in\gamma} a_{i,j}X^iY^j$. 

The polynomial $C(X,Y)$ is totally nondegenerate if the
truncation of  $C$ to any face of $\Delta(C)$ is nonsingular in
$(\mathbb C^*)^2$.
\end{defi}

A patchwork is a way to glue together the charts of real algebraic
curves.
\begin{defi}
Let $\Delta$ be a convex polygon with vertices in $\mathbb Z^2$. 
A patchwork with support
$\Delta$ is a pair $\left ( \tau, (C_i(X,Y))_{1\le i\le
  r}   \right )$ such that
\begin{itemize}
\item $\tau$ is a subdivision of $\Delta$,
\item the $(C_i(X,Y))$ are real algebraic
curves,
\item for any $2$-cell of $\tau$, there 
 is
 one and only one
  curve $C_i$ whose Newton polygon is this polygon,
\item for any edge $\gamma$ of $\tau$, the truncations on $\gamma$ of
  two 
polynomials
 whose Newton polygon contains $\gamma$ coincide.
\end{itemize}

If $\Delta$ is the triangle with vertices $(0,0)$, $(0,3)$ and
$(3n,0)$ and if the coefficient of $Y^3$ of any curve $C_i(X,Y)$ is
either $0$ or $1$,
 we say that the patchwork is trigonal of degree $n$.

\end{defi}

As we deal with nonsingular curves, we will consider patchworks
satisfying the following condition
\begin{enumerate}
\item[(1)] Each of the polynomials $C_i(X,Y)$ is totally nondegenerate.
\end{enumerate}

The next theorem is a corollary of Theorem \ref{non
  convex viro} and will be proved in section \ref{sign patchwork}.
\begin{thm}\label{viro}
Let $\Pi$ be a trigonal patchwork which satisfies hypothesis (1)
above. Then there exists a real algebraic trigonal 
curve in $\Sigma_n$ whose chart is isotopic to $\cup_{1\le i\le r}Ch(C_i)$.
\end{thm}

\vspace{2ex}
Through all this paper, we will consider trigonal patchworks which
satisfy one more condition. This condition is due to the fact that we
construct curves which have a generic position with respect to the
pencil of lines. 
This condition is a technical ingredient in the proof of Theorem
\ref{viro}. One can always perturb slightly
the curves involved in a patchwork so that they satisfy this
condition, hence it does not appear in the statement of Theorem \ref{viro}.

\begin{enumerate} 
\item[(2) ] for any $i\in\{1,\ldots,r\}$, the curves $C_i(X,Y)$ and
  $\frac{\partial 
  C_i}{\partial Y}(X,Y)$ have the maximal number
 of intersection
  points in $(\mathbb C^*)^2$ 
among curves with fixed Newton polygon $\Delta(C_i)$.
\end{enumerate}
\textbf{Remark : }The last condition simply says that the curve has
only ordinary tangency points with the vertical pencil of lines in $(\mathbb
C^*)^2$ and that none of these points lie on a  toric divisor
corresponding to 
a horizontal edge  in the toric surface associated
to $\Delta(C_i)$.

We point out a consequence of condition (2) : suppose
that $C_i(X,Y)=a_j(X)Y^k + a_{j-1}Y^{k+1} + a_{j-2}(X)Y^{k+2} \ldots
+ a_{j-l}(X)Y^{k+l}$ and that  $a_j(x_0)=0$
(resp. $a_{j-l}(x_0)=0$) with $x_0\in\mathbb R^*$, then  $a_{j-1}(x_0)\ne 0$
(resp. $a_{j-l+1}(x_0)\ne0$).

\section{The sign array of a real algebraic trigonal curve}\label{poly deg 3}

\subsection{Polynomials of degree 3 in one variable}\label{deg 3 1 var}
Let us first consider a real polynomial of degree
$3$ in one variable $C(Y)=Y^3+a_1Y^2+a_2Y+a_3$. It is well known
that the  following three numbers play an important role in the
resolution of such an equation

\begin{displaymath}
\begin{array}{ccccc}
P=a_2-\frac{a_1^2}{3},& & Q=a_3-\frac{a_2a_1}{3}+\frac{2a_1^3}{27},& & D=-4P^3-27Q^2. 
\end{array}
 \end{displaymath}

The number $D$ is the discriminant of the polynomial $C$. In particular
one has
\begin{itemize}
\item if $D>0$, then $C$ has three distinct real roots,
\item if $D<0$, then $C$ has one real root and two complex conjugated
  non-real roots,
\item if $D=0$ and $Q< 0$, then $C$ has a real double root $y_1$ and a
  simple real root $y_2$ and $y_1<y_2$,
\item if $D=0$ and $Q> 0$, then $C$ has a real double root $y_1$ and a
  simple real  root $y_2$ and $y_1>y_2$,
\item if $D=0$ and $Q=0$, then $P=0$ and $C$ has a real triple root.
\end{itemize}

\textbf{Remark : }under the change of variable $Y\mapsto -Y$,
the numbers $P$ and $D$ are invariant and the number $Q$ turns into
its opposite.

\subsection{Encoding the topology of a real algebraic trigonal
  curve}\label{topology trigonal}

Let us now consider a $\mathcal L$-nonsingular real algebraic trigonal curve $C$  in
$\Sigma_n$. Let us choose a standard coordinate system $(x,y)$ on $\Sigma_n$
such that the fiber at infinity
intersects $C$
in 3 distinct points (not necessarily real). The equation of $C$ in $\mathbb
C^2=\Sigma_n\setminus (E\cup F_\infty)$ is
$$C(X,Y)=Y^3+a_1(X)Y^2+a_2(X)Y+a_3(X)$$
where $a_j(X)$ is a polynomial of degree  $jn$. 

As in section \ref{deg 3 1 var}, define the three following polynomial in
one variable
\begin{displaymath}
\begin{array}{ccc}
P(X)=a_2(X)-\frac{a_1(X)^2}{3},&  Q(X)=a_3(X)-\frac{a_2(X)a_1(X)}{3}+\frac{2a_1(X)^3}{27},&  D(X)=-4P(X)^3-27Q(X)^2. 
\end{array}
 \end{displaymath}

The curve $C$ is $\mathcal L$-nonsingular so $D(X)$ has only simple
roots. Denote by
$x_1<x_2<\ldots<x_l$ the roots of $D(X)$. 
Denote also by $s_0$ the
sign of $D(x)$ with $x<x_1$ and $s_i$ the sign of $Q(x_i)$ for $i\in
\{1\ldots l\}$.

\begin{defi}
The sign array $[s_0,s_1\ldots s_l]$
is called the
sign array  of the curve $C$ 
in the chosen standard coordinate system.
\end{defi}

From the definition one can see that the sign array determines the
position of
a $\mathcal L$-nonsingular real algebraic trigonal curve  in
$\Sigma_n$ with respect 
to $\mathcal L$ and conversely 
provided that a standard coordinate system
has been chosen.

\begin{figure}[h]
      \centering
 \begin{tabular}{cc}
 \includegraphics[height=3.5cm, angle=0]{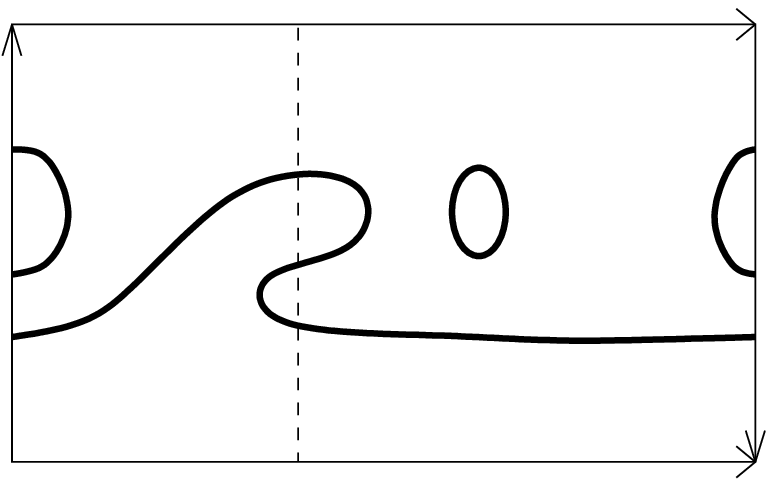}&
\includegraphics[height=3.5cm, angle=0]{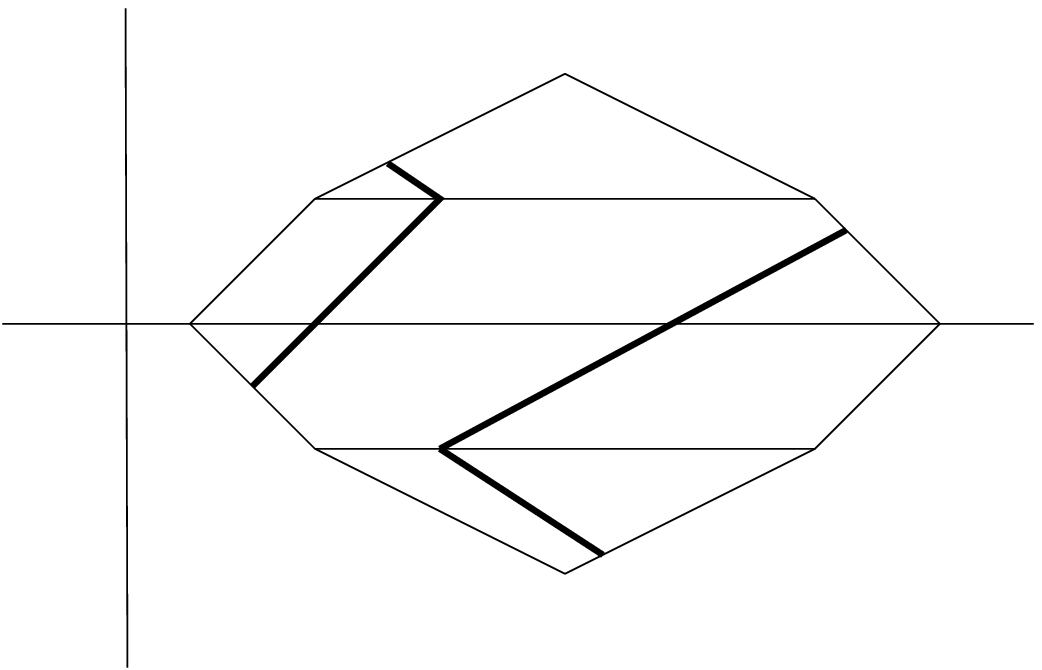}
\\a)&b)
 \end{tabular}
\caption{}
 \label{ex array}
\end{figure}

\vspace{2ex}
\textbf{Example : }
Consider the trigonal curve on $\Sigma_1$
depicted in Figure \ref{ex array}a) and suppose that a standard
coordinate system is chosen such that the vertical edges of the
rectangle represent the fiber at infinity. Then
the sign array of this trigonal curve  is $[+,+-++++]$.

\section{An order on the polygons of a subdivision}\label{pencil}
Let $\Delta$ be the triangle with vertices $(0,0)$, $(0,d)$ and
$(nd,0)$ and let $\tau$ be a subdivision of $\Delta$ into convex
integer polygons.
We 
 give an algorithm to put an order on the polygons of $\tau$. In
the Viro method, it will correspond to the order in which the polygons
are scanned by the pencil of lines. Roughly speaking, the polygons
will be scanned from the left to the right and only one polygon will
be scanned at the same time except if it contains a horizontal
edge.
If two polygons have a horizontal edge in common, they will be scanned
simultaneously by the pencil.

 \begin{figure}[h]
      \centering
 \begin{tabular}{c}
\includegraphics[height=6cm, angle=0]{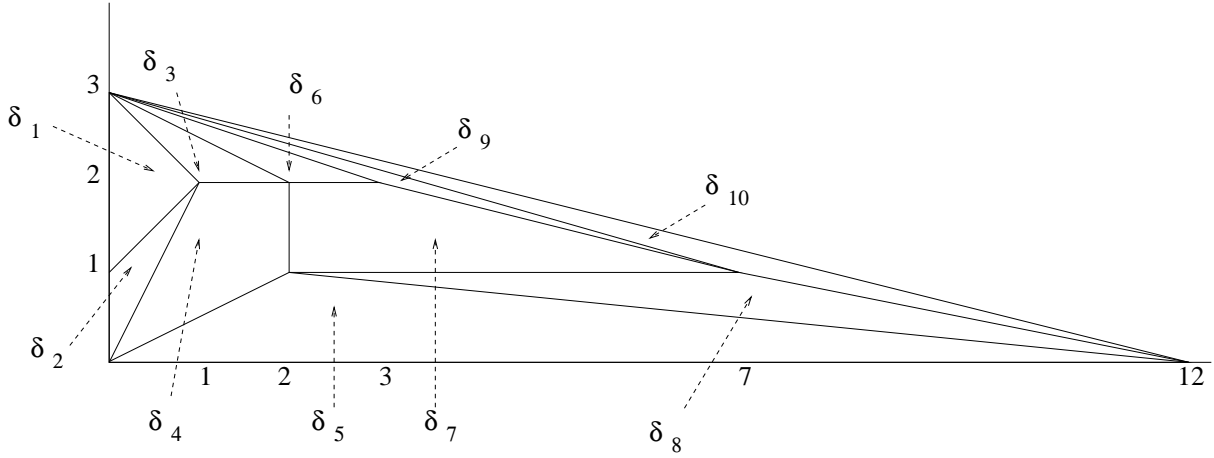}
\end{tabular}
\caption{Order on polygons}
 \label{non convex subd}
\end{figure}

\vspace{2ex}
In the procedure, $\theta$ is the current ordered subset of the
2-dimensional cells of $\tau$, and $\xi$ is the set of remaining
cells. The procedure starts with $\theta=()$ and
$\xi=\{\textrm{2-dimensional cells of }\tau\}$. Perform the following
operations until 
$\xi=\emptyset$. 

\vspace{1ex}
\textit{Step 1 : }Suppose that $\theta=(\Phi_1,\ldots,\Phi_{j-1})$ and put
\begin{itemize}
\item $U=\bigcup_{\delta\in \xi}
  \delta$,
\item $\Lambda=\{(x,y)|x=\min_{(z,y)\in U}(z)\}$ 
 
  (i.e. $\Lambda$
  is the left side of $U$),
\item $\widetilde{\xi}=\{\delta\in\xi|\delta$ has
 an edge which is contained
in $\Lambda\}$

(i.e. elements of $\widetilde{\xi}$ are the 
leftest polygons of $\xi$),
\item $\widetilde{\Lambda}=\{(\delta,(x,y))\in \widetilde{\xi}\times \Lambda | (x,y)\in\delta$, 
  $y=\max_{(z,w)\in\delta}(w)$ and $\delta$ has no horizontal edge
  containing $(x,y)\}$

 (i.e. we take 
  pairs  
$(\delta,(x,y))$
 made of a polygon $\delta$ and a distinguished vertex $(x,y)$ of $\delta$ 
such that
  $(x,y)$ is the higher point of $\delta$ and is not on one of its
  horizontal edge),
\item
  $\{(\delta_0,(x_0,y_0))\}=\{(\delta,(x,y))\in\widetilde{\Lambda}|
  y=\min_{(\delta',(z,w))\in\widetilde{\Lambda}}(w) \}$ 

(i.e. we take
  the element of $\widetilde{\Lambda}$ 
with the distinguished vertex of smaller ordinate).

\end{itemize}

\vspace{1ex}
\textit{Step 2 : }Here we look for all the polygons which will be
scanned at the same time by the pencil of lines.

Define by induction some polygons
$(\delta_{i})_{0\le i \le k}$, where $\delta_0$ has been defined
at
step~1~: if $\delta_l$ shares a horizontal edge with a polygon $\delta^\prime$, 
 distinct from $\delta_{l-1}$ if $l>0$, then put $\delta_{l+1}=
 \delta^\prime$. Else, put $k=l$ and $\Phi_j=\{\delta_i, 0\le i \le
      k\}$.

\vspace{1ex}
\textit{Step 3 : }Put $\xi=\xi\setminus \Phi_j$. If 
$x_0>1$ then
also put
$\theta=(\Phi_1,\ldots,\Phi_{j})$.

\vspace{2ex}
\textbf{Remark : }in this algorithm, we 
``throw away''  polygons
corresponding to curves of degree 1 in $Y$ and which
 do not have a horizontal edge of height 1.

\vspace{2ex}
\textbf{Example : }Applying this algorithm to the
(non-convex) 
subdivision shown in
Figure~\ref{non convex subd}, one obtains 
$$\theta=(\{\delta_1\},
\{\delta_2\}, \{\delta_3, \delta_4\}, 
\{\delta_6,\delta_7,\delta_8\ \}, \{\delta_9\}, \{\delta_{10}\}).$$

 \begin{figure}[h]
      \centering
 \begin{tabular}{c}
\includegraphics[height=4cm, angle=0]{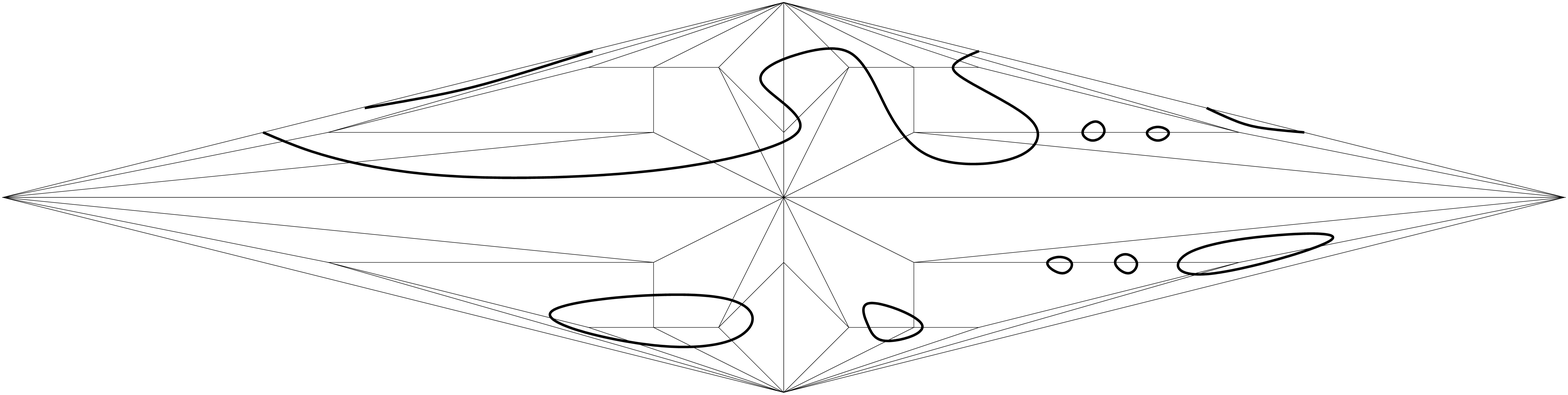}
\end{tabular}
\caption{}
 \label{non convex patch}
\end{figure}

Let $\Phi$ be an element of $\theta$ and $\delta$ a polygon of $\Phi$.
\begin{defi}
  The height of $\delta$ is defined as $\max_{(x,y)\in\delta}(y)$.
  The height of $\Phi$ is the maximum of the heights of the polygons of
  $\Phi$.
\end{defi}
Note that by construction, all the elements of $\Phi$ 
have height at least 2.

\section{The sign array associated to a trigonal patchwork}\label{sign
  patchwork}

To a trigonal  patchwork, we associate a sign array encoding the position
that 
the curve we want to construct 
would have
with respect to $\mathcal L$.
This array is constructed studying where should be located
the roots of the discriminant of the desired curve, and
what should be the sign of $Q(X)$ at those points.

\begin{figure}[h]
      \centering
 \begin{tabular}{ccc}
\includegraphics[height=6cm, angle=0]{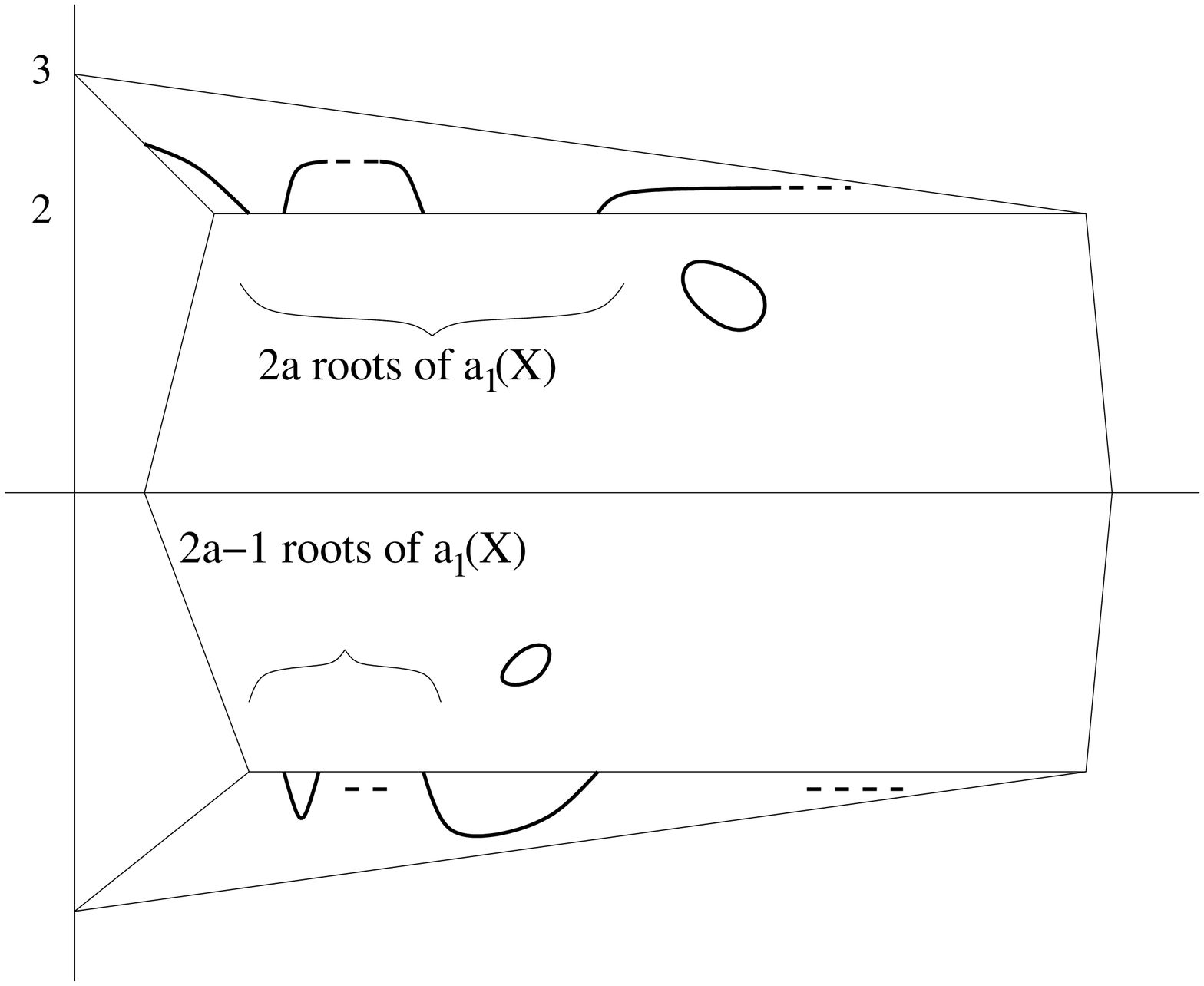}&
\hspace{5ex}&

\includegraphics[height=6cm, angle=0]{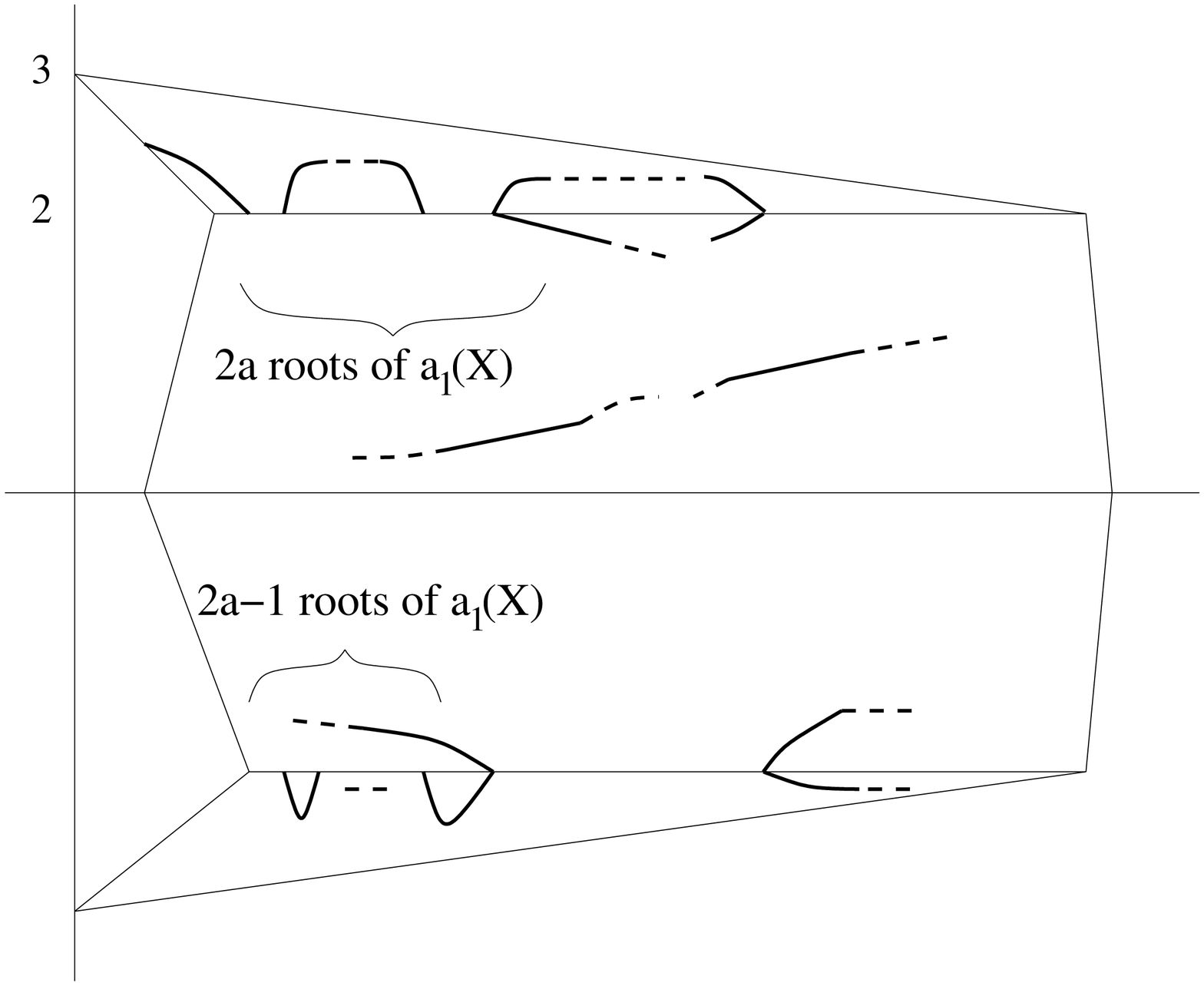}

\\a)&&b)

\end{tabular}
\caption{}
 \label{horiz 2}
\end{figure}

Let $\Pi$ be a trigonal patchwork satisfying conditions (1) and (2) of
section \ref{state viro}. 
Following what happens when the subdivision is convex (see \cite{Lop}),
the roots $x_0\in\mathbb R$ of $D(X)$ 
 will be obtained out of $\Pi$ in exactly three different ways:

\begin{itemize}
\item any root in $\mathbb R^*$ of the discriminant a curve $C_i$ of
  the patchwork gives rise to  such a root
 $x_0$,

\item suppose there
is a horizontal edge of height 1 
 contained in a polygon $\delta_1$  of
  height at least 2 and in a polygon $\delta_2$ of height 1. The
  polygon 
  $\delta_1$ corresponds in the patchwork to a curve of equation
  $\alpha Y^3+a_1(x)Y^2+a_2(X)Y$ and the polygon
  $\delta_2$ corresponds  to a curve of equation
  $a_2(X)Y+a_3(X)$.

Then a
 root of $x_1\in\mathbb R^*$ of $a_2(X)$ such that $a_1(x_1)$
  and $a_3(x_1)$ have the same sign (see Figure \ref{ex array}b))
gives rise to
 two such roots
 $x_0$.

\item suppose there is a  horizontal 
edge of height 2 
  contained in a polygon $\delta_1$  of
  height 3 and in a polygon $\delta_2$ of height 2. The
  polygon 
  $\delta_1$ corresponds in the patchwork to a curve of equation
  $Y^3+a_1(x)Y^2$ and the polygon
  $\delta_2$ corresponds  to a curve of equation
  $a_1(X)Y^2+a_2(X)Y+a_3(X)$.

Then a root
$x_1\in\mathbb R^*$ of $a_1(X)$ such that $a_2(x_1)$
  is positive (see Figure \ref{horiz 2}a))
gives rise to
 two such roots
 $x_0$.
\end{itemize}

First, we associate two lists of signs to each element of $\theta$.

\subsection{Lists of signs associated to an element of 
$\theta$}
Let $\Phi_i$ be the $i^{th}$ element of
 $\theta$ and let be $a_1(X)$,
$a_2(X)$ and $a_3(X)$ all the polynomials that occur in the polygons
of $\Phi_i$. The polynomial $a_j(X)$ is the coefficient of $Y^{(3-j)}$
for some polygon of $\Phi_i$.

We define $x^-_{i,1}<\ldots< x^-_{i,u_i}$ (resp. $x^+_{i,1}<\ldots<
x^+_{i,v_i}$) to be all the numbers in $\mathbb R^*_-$ (resp. $\mathbb
R^*_+$) which verify one of the following conditions
\begin{itemize}
\item $x^s_{i,l}$ is a root of the discriminant of the curve
  corresponding to one of the polygons of $\Phi_i$,

\item two polygons of $\Phi_i$ have a common horizontal edge of height
  1, and $x^s_{i,l}$ is a root of $a_2(X)$ such that $a_1(x^s_{i,l})$
  and $a_3(x^s_{i,l})$ have the same sign (see Figure \ref{ex array}b)).

\item two polygons of $\Phi_i$ have a common horizontal edge of height
  2, and $x^s_{i,l}$ is a root of $a_1(X)$ such that $a_2(x^s_{i,l})$
  is positive  (see Figure \ref{horiz 2}).
\end{itemize}

\subsubsection{$\Phi_i$ is of height 3 and does not contain any
  polygon of height 2 }

That means that $\Phi_i$ contains a polygon corresponding to a curve
of degree $3$ and that
no one of the $x^s_{i,j}$ comes from
an edge of height 2.
Denote by $\sigma^s_i$ the sign of
  $Q_i(x^s_{i,l})$ and perform the following substitutions in the
lists $x^-_{i,1}x^-_{i,2}\ldots x^-_{i,u_i}$  and $x^+_{i,1}x^+_{i,2}\ldots x^+_{i,u_i}$
\begin{itemize}
\item if $x^s_{i,l}$ is a root of the discriminant of the curve of
  degree 3, then replace $x^s_{i,l}$ by $\sigma^s_l$,

\item if $x^s_{i,l}$ is a root of $a_2(X)$, then
  replace $x^s_{i,l}$ by $\sigma^s_l\sigma^s_l$ (see Figure \ref{ex array}b)).
\end{itemize}

\vspace{2ex}
\textbf{Example : }The two lists of signs  associated to $\{\delta_1\}$ in the
patchwork depicted 
in Figure \ref{non convex patch} are $-+$ and the empty list.

\subsubsection{$\Phi_i$ contains a
  polygon of height 2 }

That means that $\Phi_i$ contains a polygon corresponding to a curve
of degree $2$.
Denote by $\sigma^s_i$  the sign of
  $a_1(x)$ for  $x\in\mathbb R^*_s$ small enough, and by  $a^s_l$
the number of roots of $a_1(X)$  counted with multiplicity which are
strictly between $0$ and $x^s_{i,l}$.
Define $\sigma'^s_l=(-1)^{a^s_l}\sigma^s_i$ and $\sigma''^s_l=-\sigma'^s_l$ and perform the
following substitutions in the
lists $x^-_{i,1}x^-_{i,2}\ldots x^-_{i,u_i}$  and
$x^+_{i,1}x^+_{i,2}\ldots x^+_{i,u_i}$ 
\begin{itemize}
\item if $x^s_{i,l}$ is a root of the discriminant of the curve of
  degree 2, then  replace $x^s_{i,l}$ by $\sigma'^s_l$ (see Figure \ref{horiz 2}a)) 

\item if $x^s_{i,l}$ is a root of $a_2(X)$, then
  replace $x^s_{i,l}$ by $\sigma'^s_l\sigma'^s_l$, 

\item if $x^s_{i,l}$ is a root of $a_1(X)$, then
  replace $x^s_{i,l}$ by $\sigma''^s_l\sigma'^s_l $ (see Figure \ref{horiz 2}b)).

\end{itemize}

\vspace{2ex}
\textbf{Example : }The two lists of signs  associated to
$\{\delta_6,\delta_7,\delta_8\}$ in the 
patchwork depicted 
in Figure \ref{non convex patch} are the empty list and $-+----------$

\subsection{The final sign array}

So for each element $\Phi_i$ of $\tau$, we have defined two lists of signs
$\sigma^-_{i,1}\ldots \sigma^-_{i,u'_i}$ and $\sigma^+_{i,1}\ldots
\sigma^+_{i,v'_i}$ .  
 
Define $b_{(3-i)n,i}$ to be the coefficient of the monomial
$X^{(3-i)n}Y^i$ in the patchwork and put
$\widetilde C(Y)=\sum_{i=0}^3b_{(3-i)n,i}Y^i$. Denote by $s_0$ the
sign of the discriminant of $\widetilde C(Y)$.

\begin{defi}
The sign array associated to the patchwork is
$$[s_0,\sigma^-_{k,1}\sigma^-_{k,2}\ldots
  \sigma^-_{k,u'_k}\sigma^-_{k-1,1} \ldots
  \sigma^-_{k-1,u'_{k-1}}\ldots \ldots\sigma^-_{1,1}\ldots \sigma^-_{1,u'_1}
\sigma^+_{1,1}\ldots
\sigma^+_{1,v'_1}\sigma^+_{2,1}\ldots
\sigma^+_{2,v'_2}\ldots\ldots\sigma^+_{k,1}\ldots
\sigma^+_{k,v'_k}]
 $$
\end{defi}

\vspace{2ex}
\textbf{Example : }The sign array associated to the trigonal patchwork depicted
in Figure \ref{non convex patch} is 
$$[+,+--+---+------------].$$

\vspace{2ex}
We can now state our main theorem 
whose proof will be given in section~\ref{proof}, Proposition~\ref{non convex viro prop}.
\begin{thm}\label{non convex viro}
Let $\Pi$ be a trigonal patchwork of degree $n$ which satisfy the hypothesis (1) and (2) of section
\ref{state viro}. Then there exists a $\mathcal 
L$-nonsingular real algebraic trigonal curve in $\Sigma_n$ realizing the sign array
associated to $\Pi$.
\end{thm}

\vspace{2ex}
Using Theorem \ref{non convex viro}, we can now prove Theorem
\ref{viro}.
\vspace{1ex}

\textit{Proof of Theorem \ref{viro}: }If the hypothesis (1) and (2) of
section \ref{state viro} are fulfilled, then Theorem \ref{viro} is a
consequence of Theorem \ref{non convex viro}. Indeed, any $\mathcal
L$-nonsingular real algebraic trigonal curve realizing the sign array
associated to $\Pi$ has a chart which is isotopic to the union
of the charts used in the patchwork.

If only the hypothesis (1)  is fulfilled, one can perturb slightly the
coefficients of the curves used in the patchwork in order that the
hypothesis (2) of
section \ref{state viro} is fulfilled. Then apply Theorem \ref{non
  convex viro} to this patchwork. Perturbing in different ways the
initial patchwork, one can construct $\mathcal
L$-nonsingular real algebraic trigonal curve realizing different
positions with respect to the pencil $\mathcal L$. However, all these
curves have isotopic charts 
and this isotopy type is
the one required by Theorem \ref{viro}\findemo

\section{Signed real rational graphs associated to trigonal curves and to curves of bidegree $(2,n)$}\label{graph}

\subsection{Signed real rational graphs}

Here we recall some fact about real rational graphs. We refer to
\cite{Br2}, \cite{Br1}, \cite{O3} and \cite{Shap1} for more details
and proofs.

\begin{figure}[h]
      \centering
 \begin{tabular}{ccc}
\includegraphics[height=3.5cm, angle=0]{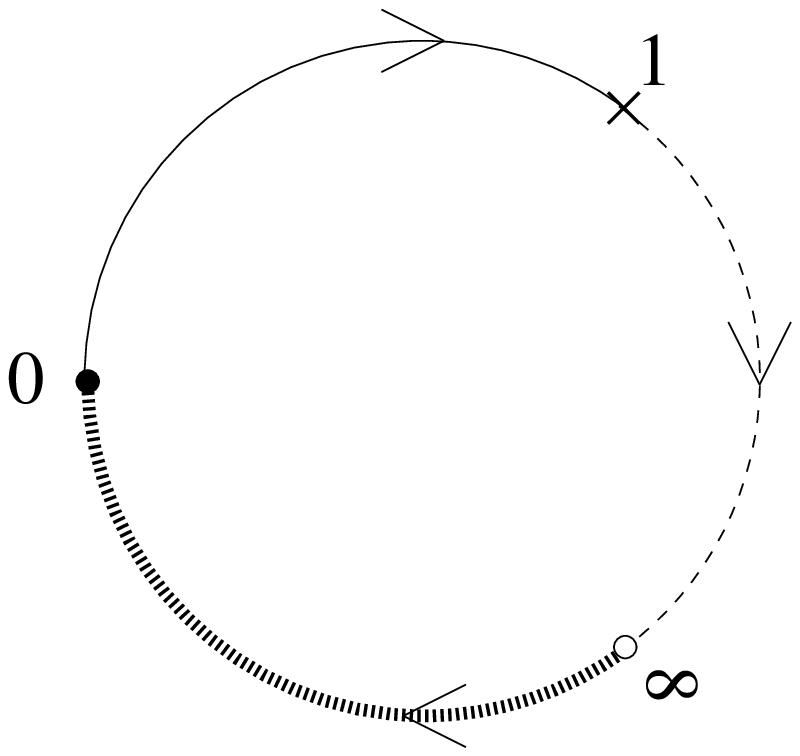}&
\hspace{15ex}  &
\includegraphics[height=3.5cm, angle=0]{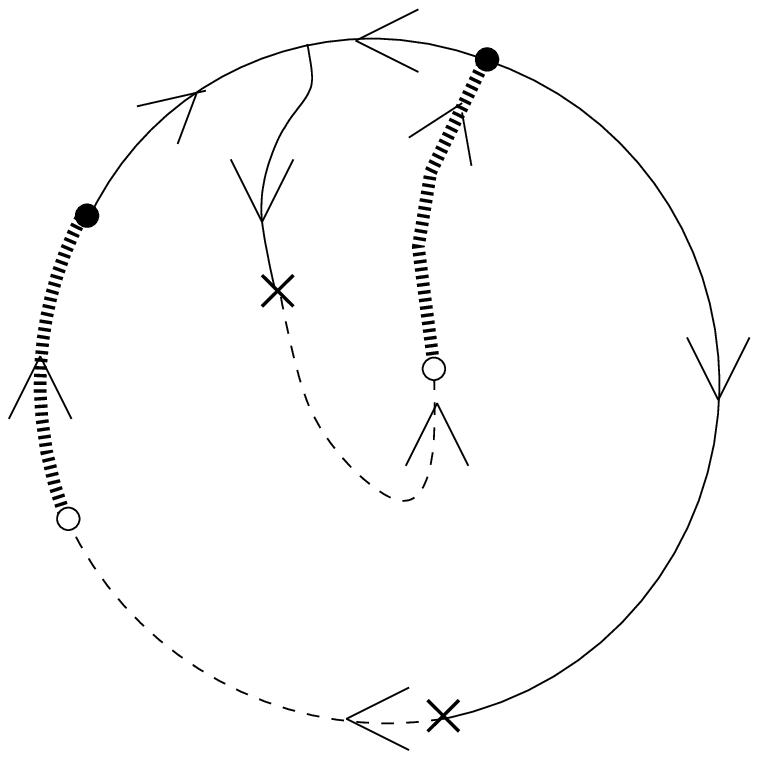}
\\ a)&&b)
 \end{tabular}
\caption{}
 \label{ex array2}
\end{figure}

Color and orient
$\mathbb RP^1$ as depicted in Figure \ref{ex array2}a).
\begin{defi}\label{defi rat grap}
Let $\Gamma$ be a graph
on $\mathbb CP^1$ invariant under the action
of the complex conjugation and $\pi:\Gamma\to\mathbb RP^1$ a
continuous map.
 Then the coloring and
orientation of $\mathbb RP^1$ shown in Figure \ref{ex array2}a) defines
a coloring and an orientation of $\Gamma$ via $\pi$.

\noindent The graph $\Gamma$ equipped with this coloring and this
orientation is called a real rational graph if 
\begin{itemize}
\item any vertex of $\Gamma$ has an even valence,
\item for any connected
component $W$ of $\mathbb CP^1\setminus\Gamma$, the map
$\pi_{|\partial W}$ is a covering of $\mathbb RP^1$ of degree $d_W$, 
\item  for any connected
component $W$ of $\mathbb CP^1\setminus\Gamma$, the orientation
induced by $\pi$ on $\partial W$ is also induced by an orientation of $W$.
\end{itemize} 
The sum of the degrees $d_W$ for all connected
component $W$ of $\{Im(z)>0\}\setminus\Gamma$ of is called the
degree of $\Gamma$. 
\end{defi}

Since a real rational graph is invariant under the complex conjugation,
we will draw in this article only one of the two halves of  real
rational graphs. This half will be drawn in a disk and the boundary of
this disk will be $\mathbb RP^1$.

Let $f:\mathbb CP^1\to\mathbb CP^1$ be a
real rational map
 of degree $d$, and let $\Gamma$ be $f^{-1}(\mathbb RP^1)$ with   the
coloring and the orientation induced by those chosen on $\mathbb
RP^1$. Then $\Gamma$ is a real rational graph of degree $d$.

\begin{defi}
One says that $\Gamma$ is the real rational graph
associated to $f$.
\end{defi}

\textbf{Example : }the real rational graph associated to the real rational
map $f(X,Z)=\frac{(X-Z)(X-2Z)^2}{(X+3Z)(X^2+XZ+Z^2)}$ is depicted in
Figure \ref{ex array2}b).

\vspace{2ex}
The next theorem shows that the converse is also true : one can
construct real rational map out of a real rational graph.

\begin{thm}\label{constr rat}
Let  $\Gamma$ be a real rational graph of degree $d$. Then there
exists a 
real  
rational map
of degree $d$
 whose real rational graph is equivariantly
isotopic to $\Gamma$.

Such a real rational map is said to realize $\Gamma$.\findemo 
\end{thm}

One of the  steps in the proof of our main Theorem is the gluing
of some real rational graphs. These 
gluings will not depend only on the real rational graphs, but also on
some signs of polynomials which define it and on some coordinate
system on $\mathbb CP^1$. This motivates the notions of  marked and
signed real 
rational graphs.

\begin{defi}
A marked real rational graph is a real rational graph on $\mathbb
CP^1$ where some coordinate system is fixed. In particular, the points 
 $0$ and $\infty$ are prescribed on such a graph.

\vspace{2ex}
Let $\Gamma$ be a marked real rational graph and $f$ a real rational
map which realizes $\Gamma$.
The marked graph $\Gamma$ is called a signed real rational graph if
each connected component of $(\Gamma\cap \mathbb RP^1)\setminus
(f^{-1}(\{0,1,\infty\})\cup \{0,\infty\})$ is
enhanced with a pair of signs.
\end{defi}

As explained above, we will represent a real rational graph by a graph
in a disk lying in $\mathbb R^2$. A marked real rational graph will be
represented in a disk which is symmetric with respect to the
$x$-axis.
If $(x_1,0)$ and $(x_2,0)$ are the two intersection
points of $\{y=0\}$ and the circle bounding this disk with $x_1<x_2$,
then $(x_1,0)$ (resp. $(x_2,0)$) will be $0$ (resp. $\infty$) on  the
marked graph. The positive (resp. negative) points of $\mathbb RP^1$ will
be drawn on the half plane $\{(x,y)|y>0\}$ (resp. $\{(x,y)|y<0\}$) and
is called the positive (resp. negative) part of the graph.

\vspace{2ex}
In the rest of this paper, we will not draw the arrows on the real
rational graphs.

\vspace{2ex}
\textbf{Example : }A signed real rational graph is depicted on figure
\ref{ex rat graph}b).

\subsection{Signed real rational graph associated to a real algebraic
  trigonal curve in $\Sigma_n$}
Let $C$ be a real algebraic trigonal curve in $\Sigma_n$. Then $C$ has
the following equation in $\Sigma_n\setminus E$
$$C(X,Y,Z)=Y^3+a_1(X,Z)Y^2+a_2(X,Z)Y+a_3(X,Z)$$
where $a_j(X,Z)$ is a homogeneous polynomial of degree $jn$. To such a
curve, we associate the real rational graph
 $\Gamma(C)$
 given by the map
\begin{displaymath}
\begin{array}{ccccc}
f & : & \mathbb CP^1 &  \longrightarrow & \mathbb CP^1\\
& & [X:Z]& \longmapsto &\frac{-4P(X,Z)^3}{27Q(X,Z)^2}
\end{array}
\end{displaymath} 
where the polynomials $P$ and $Q$ are those defined in section
\ref{topology trigonal}.
\begin{defi}
The graph $\Gamma(C)$ is called the real rational graph associated to
$C$.
\end{defi}

The preimages of $1$ by $f$ correspond to the roots of $D(X,Z)$, the
discriminant of $C$.
So if the curve $C$ is $\mathcal L$-nonsingular, according to section \ref{poly deg
  3}, one can recover the position of $C$ in $\Sigma_n$ with respect
to $\mathcal L$ up to the transformation $Y\mapsto -Y$ only from
$\Gamma\cap\mathbb RP^1$.

 \begin{figure}[h]
      \centering
 \begin{tabular}{ccc}
\includegraphics[height=4.5cm, angle=0]{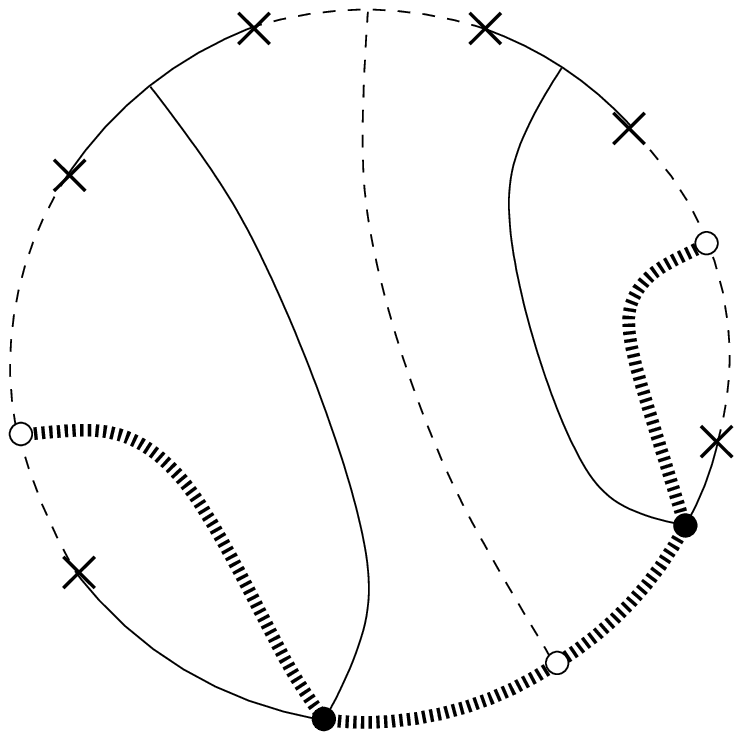}&&
\includegraphics[height=4.5cm, angle=0]{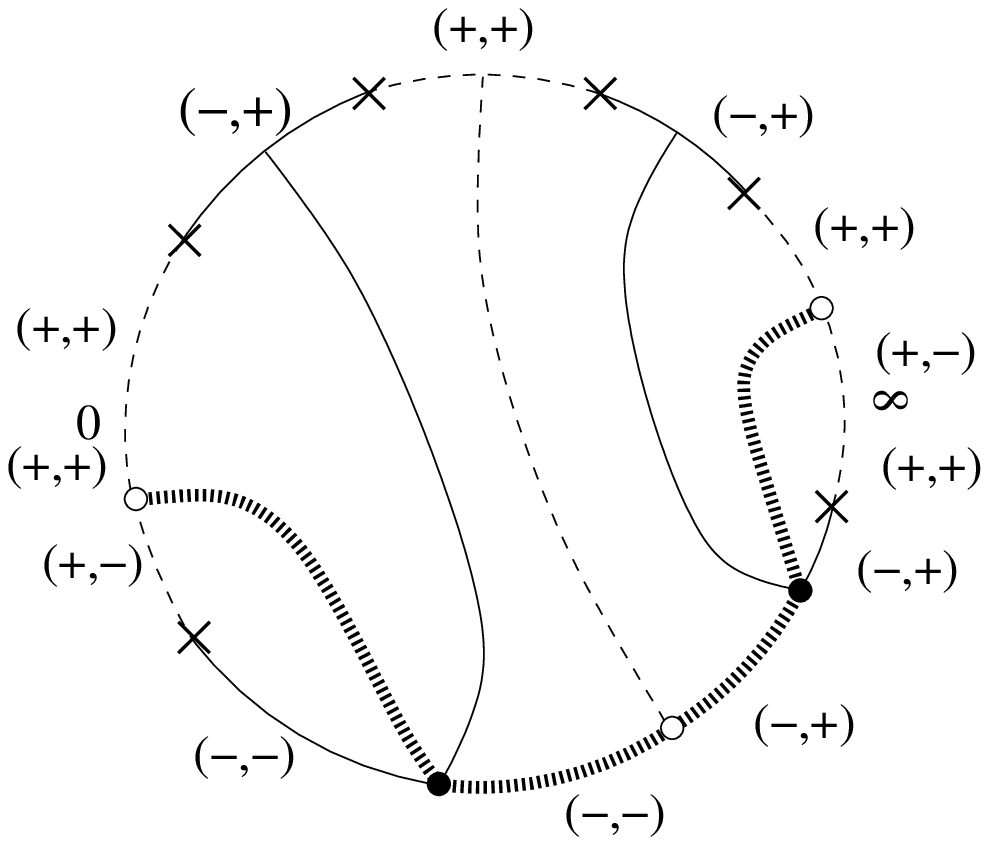}

\\ a)&&b)
\end{tabular}
\caption{}
 \label{ex rat graph}
\end{figure}

\vspace{2ex}
\textbf{Example : }the real rational graph associated to the trigonal
curve in $\Sigma_1$ depicted on Figure \ref{ex array}a) is depicted in
Figure \ref{ex rat graph}a).

\vspace{2ex}
Conversely, suppose that using Theorem \ref{constr rat}, we have
constructed a real rational map 
$f(X,Z)=\frac{-4P(X,Z)^3}{27Q(X,Z)^2}$ out of 
a real rational graph $\Gamma$. Then, $\Gamma$ is the real rational
graph associated to the trigonal curve
$Y^3+P(X,Z)Y+Q(X,Z)$.

\vspace{2ex}
Up to rotation, this real rational graph does not depend on the
coordinate system we choose on $\Sigma_n\setminus E$. If now we take
into account the coordinate system, we obtain a marked real rational
graph. Let us 
 turn it into 
 a signed real rational graph 
$\Gamma_\pm(C)$~:
 the pair of signs on  each connected component of
$(\Gamma\cap \mathbb RP^1)\setminus 
(f^{-1}(\{0,1,\infty\})\cup \{0,\infty\})$  is the
pair formed by the sign of $D(X)$ and the sign of $Q(X)$ on this
component. 
\begin{defi}
The graph $\Gamma_\pm(C)$ 
is called the signed real rational graph associated to the
polynomial $C$ .
\end{defi}

\vspace{2ex}
\textbf{Example : }Suppose that we have fixed the coordinate system in
$\Sigma_1$ such that the vertical edges of the rectangle in Figure
\ref{ex array}a) are the fiber at infinity and the fiber in dotted
line is the fiber at $0$. Then  $\Gamma_\pm(C)$ is depicted in Figure \ref{ex rat graph}b).

\vspace{2ex}
Note that the signed real rational graph associated to a trigonal
curve in some coordinate system 
can be extracted from the knowledge of the real rational graph
associated to this curve and the topology of the curve.
For this reason, we won't write explicitly the pairs of signs on
signed real rational graph associated to trigonal curves in the rest
of this paper.

\vspace{2ex}
Here we give a list of necessary and sufficient conditions for a signed
real rational graph to be associated to some $\mathcal L$-nonsingular
real algebraic trigonal curve.
The proof can be found in \cite{O3}.

\begin{prop}\label{iff trigo}
Let $\Gamma_{\pm}$ be a signed real rational graph and let 
$\pi:\Gamma_{\pm}\to\mathbb RP^1$ be a continuous map as in definition
  \ref{defi rat grap}. 
Then there exists a
$\mathcal L$-nonsingular 
real algebraic trigonal curve in $\Sigma_n$ such that
$\Gamma_{\pm}(C)$ is equivariantly isotopic to $\Gamma_{\pm}$ if and
only if the following conditions are fulfilled
\begin{itemize}
\item $\Gamma_{\pm}$ is of degree $6n$,
\item any preimage of $0$ has an order which is divisible by $3$,
\item any preimage of $\infty$ has an order which is divisible by $2$,
\item any preimage of $1$ is of order $1$,
\item any pair of sign $(+,s)$ labels a connected
  component of $\pi^{-1}(]\infty,0[)$ or $\pi^{-1}(]\infty,0[)$, 
\item any pair of sign $(-,s)$ labels a connected
  component of $\pi^{-1}(]1,\infty[)$,
\item when passing through a preimage of $\infty$ of order $2a$, the
  pair of sign $(s',s)$ becomes $(s',(-1)^as)$.~\findemo
\end{itemize}

\end{prop}

\subsection{Signed real rational graph associated to a real algebraic
   curve of bidegree (2,n) in $\Sigma_n$}\label{rat bigon}

 \begin{figure}[h]
      \centering
 \begin{tabular}{ccc}
 \includegraphics[height=3.5cm, angle=0]{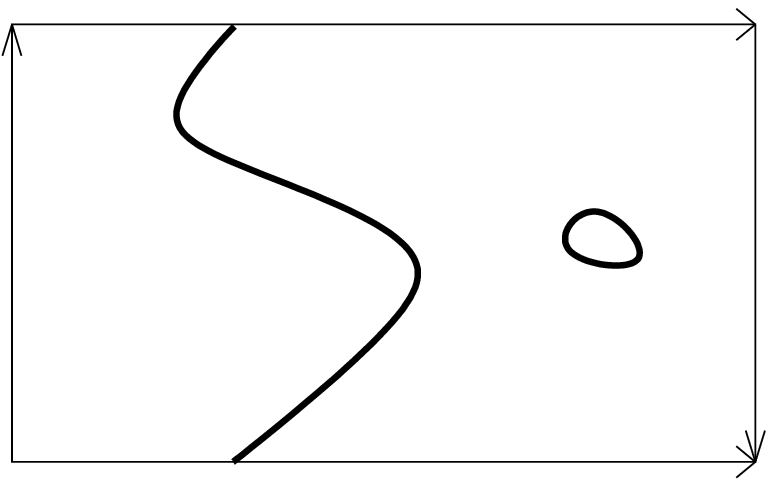} &&
\includegraphics[height=3.5cm, angle=0]{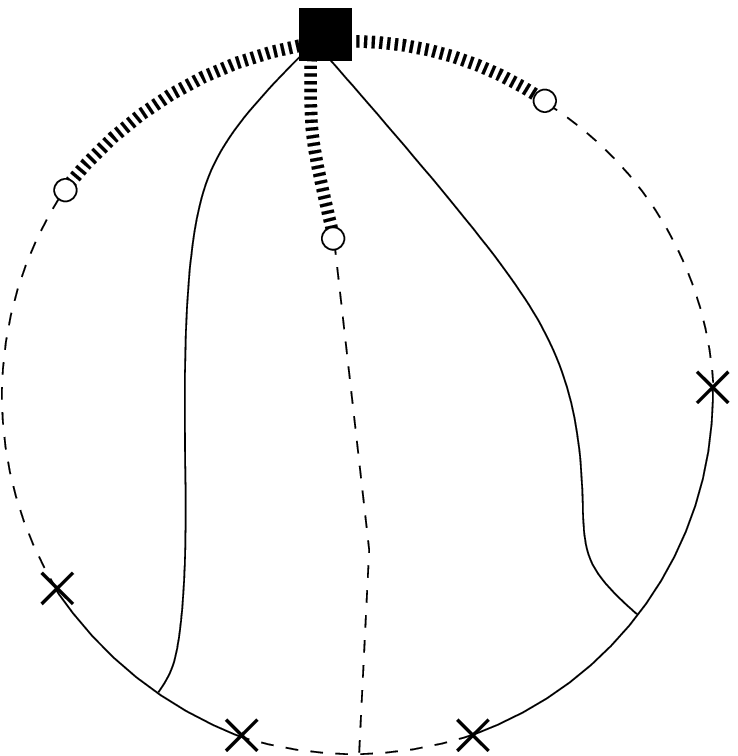}

\\ a)&&b)
\end{tabular}
\caption{}
 \label{ex rat graph 2}
\end{figure}

Let $C$ be a real algebraic  curve of bidegree $(2,n)$ in
$\Sigma_n$. So $C$ has 
the following equation in $\Sigma_n\setminus E$
$$C(X,Y,Z)=a_1(X,Z)Y^2+a_2(X,Z)Y+a_3(X,Z)$$
where $a_j(X,Z)$ is a homogeneous polynomial of degree $jn$. The
intersections of the curve $C(X,Y,Z)$  with the exceptional divisor correspond
to the roots of $a_1(X,Z)$.
To such a
curve, we associate the real rational graph $\Gamma(C)$ associated to the map
\begin{displaymath}
\begin{array}{ccccc}
f & : & \mathbb CP^1 &  \longrightarrow & \mathbb CP^1\\
& & [X:Z]& \longmapsto &\frac{-a_1(X,Z)^4}{a_2(X,Z)^2-4a_1(X,Z)a_3(X,Z)-a_1(X,Z)^4}
\end{array}.
\end{displaymath} 
\begin{defi}
The graph $\Gamma(C)$ is called the real rational graph associated to
$C$.
\end{defi}
The real rational graph associated to the curve $C$ appears as a
 limit of  the real
rational graph associated to a real trigonal algebraic curve when its
coefficient of $Y^3$ tends to $0$.

The preimages of $1$ by $f$ correspond to the roots of $a_2(X,Z)^2-4a_1(X,Z)a_3(X,Z)$, the
discriminant of $C$.

\vspace{2ex}
\textbf{Example : }the real rational graph associated to the 
curve of bidegree $(1,2)$ in $\Sigma_1$ depicted in Figure \ref{ex rat
  graph 2}a) is depicted in
Figure \ref{ex rat graph 2}b).

\vspace{2ex}
As in the case of trigonal curves,  the real rational graph associated
to $C$ does not depend on the
coordinate system we choose on $\Sigma_n\setminus E$,up to
rotation. Taking into account the coordinate system chosen to write
the equation of $C$, we obtain a marked real rational graph. Let us
turn it into
 a signed real rational graph $\Gamma_\pm(C)$~: the pair of signs on  each connected component of
$(\Gamma\cap \mathbb RP^1)\setminus 
(f^{-1}(\{0,1,\infty\})\cup \{0,\infty\})$  is the
pair formed by the sign of $a_2(X,Z)^2-4a_1(X,Z)a_3(X,Z)$ and the sign
of $  a_1(X)$ on
this 
component.

\begin{defi}
The graph $\Gamma_\pm(C)$ 
is called the signed real rational graph associated to the
polynomial $C$ .
\end{defi}

This definition of $\Gamma_\pm$ for a curve of bidegree $(2,n)$ is
motivated by the fact that when we 
will glue together the curves of the patchwork in
section \ref{gluing graphs}, the topology coming from a curve of
degree 2
will be
 encoded by the sign of $a_1(X)$ (see Figure \ref{horiz 2}).

\vspace{2ex}
Note that the knowledge of the real rational graph and the topology of
$C$ does not allow one to recover the sign real rational
graph. Indeed, the sign of $a_1(X)$ can be recovered neither from
the real rational graph nor from the topology of the curve.

However, when the curve $C$ will be 
 used
 in a patchwork, the sign of
some other polynomials will 
 allow
 us to recover the sign of $a_1(X)$.

\section{Associating a signed real rational graph to any element of
  $\theta$}\label{constr graph}

 \begin{figure}[h]
      \centering
 \begin{tabular}{ccc}
\includegraphics[height=3.5cm, angle=0]{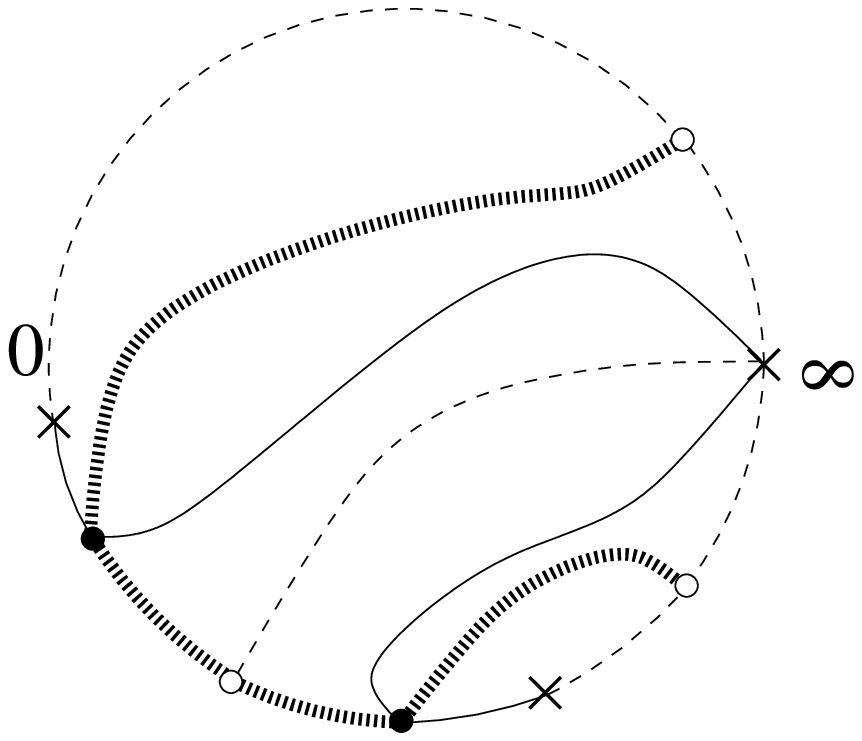}&

\includegraphics[height=3.5cm, angle=0]{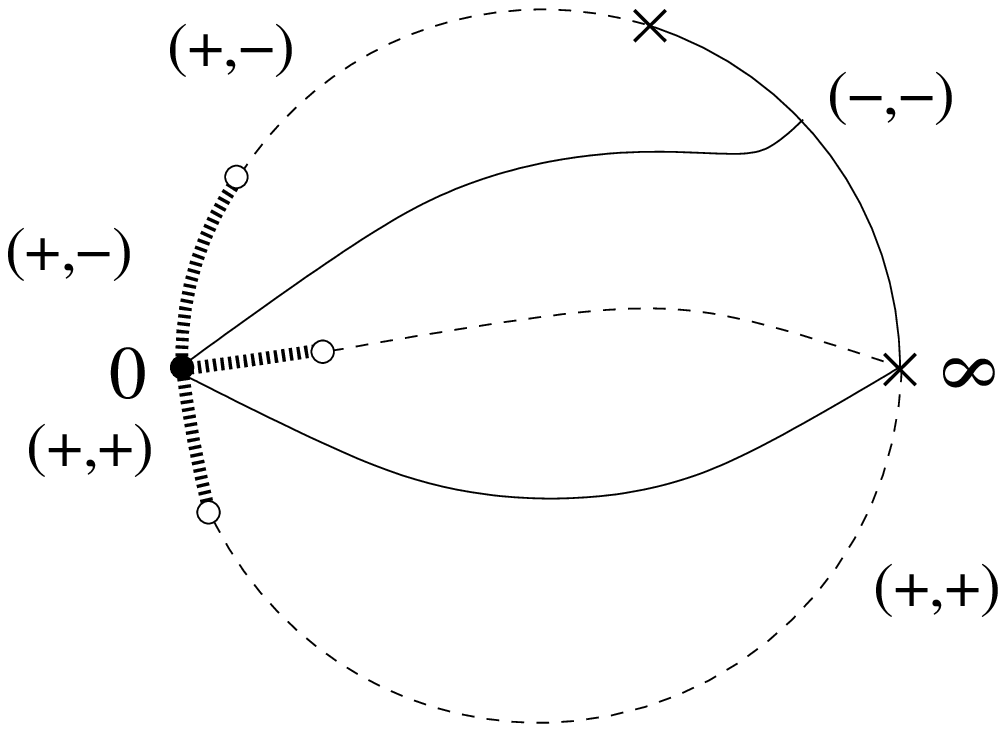}&
\includegraphics[height=3.5cm, angle=0]{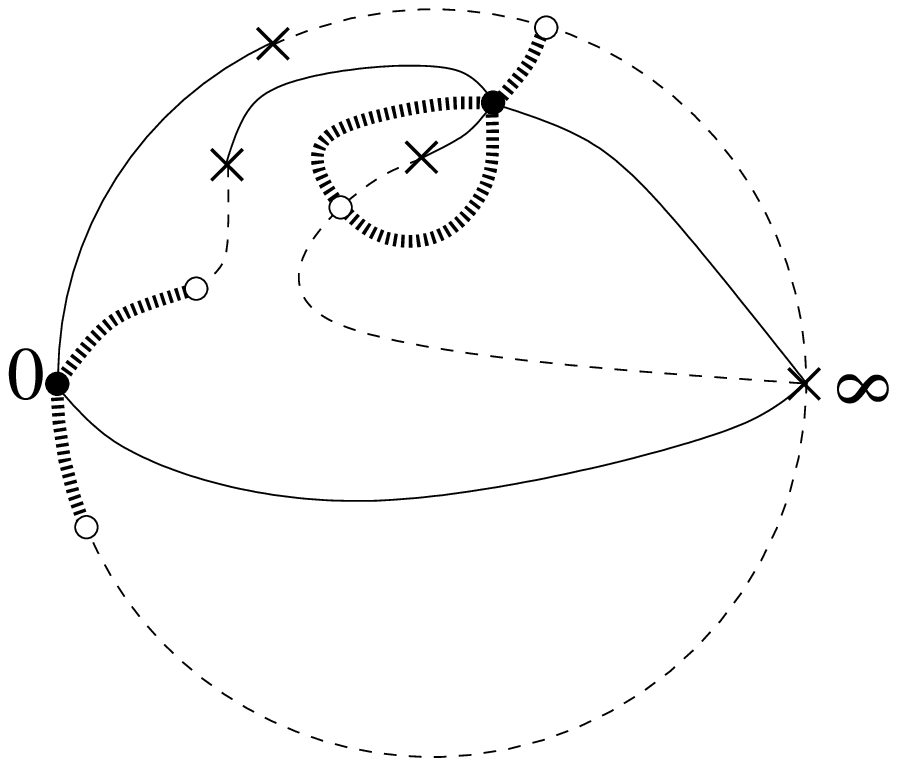}

\\ a)$\{\delta_1 \}$ &b)$\{\delta_2 \}$ &c)$\{\delta_3, \delta_4\}$
\end{tabular}
\caption{}
 \label{real graph}
\end{figure}

 \begin{figure}[h]
      \centering
 \begin{tabular}{ccc}
\includegraphics[height=3.5cm, angle=0]{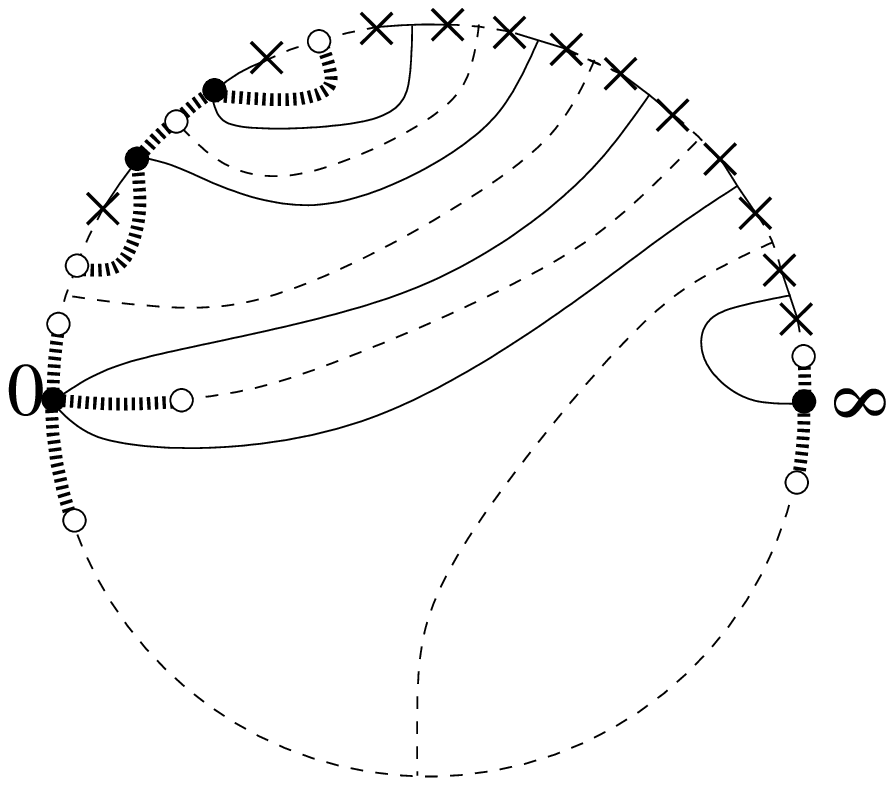}&

\includegraphics[height=3.5cm, angle=0]{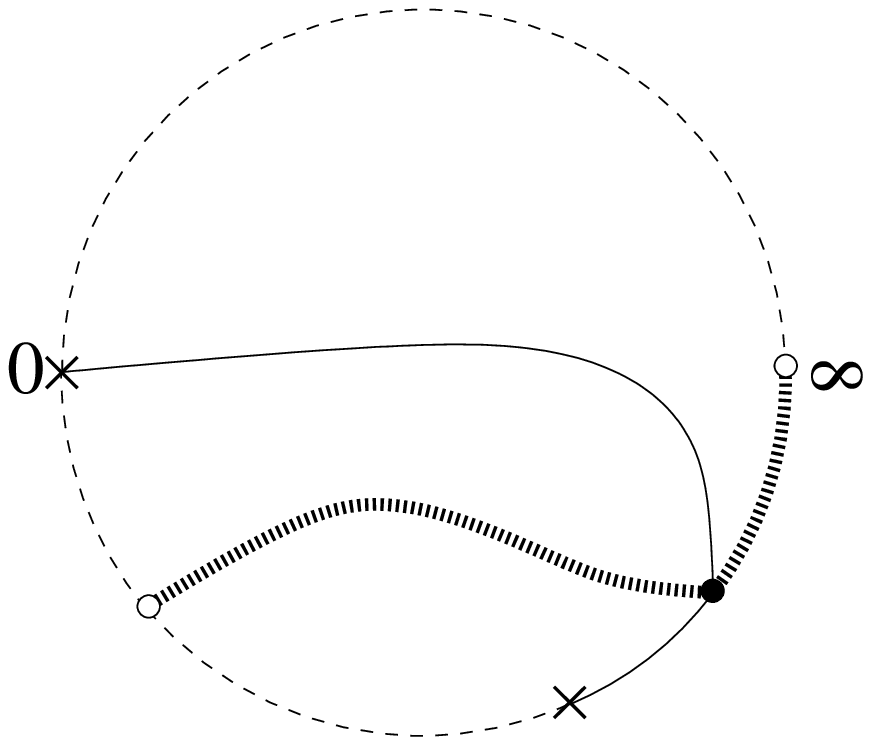}&

\includegraphics[height=3.5cm, angle=0]{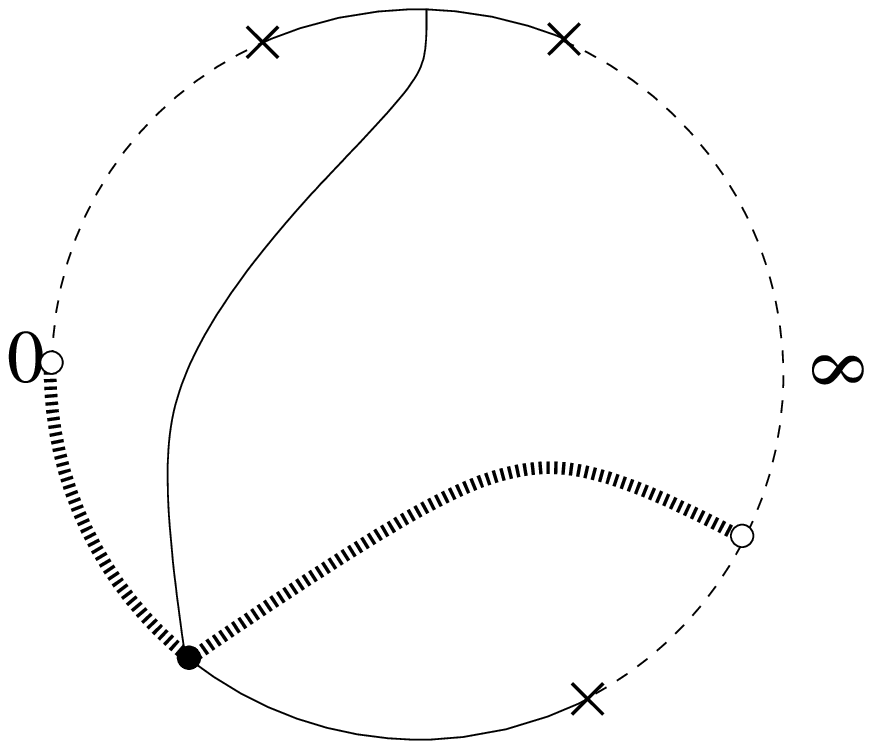}
\\ a)$\{\delta_6,\delta_7,\delta_8\}$ &b)$\{\delta_9 \}$ &c)$\{\delta_{10} \}$
\end{tabular}
\caption{}
 \label{real graph 2}
\end{figure}

Let $\Pi$ be a trigonal patchwork of degree $n$ satisfying conditions
(1) and (2) of section \ref{state viro} and let
$\theta=(\Phi_1,\ldots,\Phi_k)$ be the $k$-tuple constructed in
section \ref{pencil}. 
We are going to associate a signed real rational graph
$\Gamma_\pm(\Phi_i)$ to
 each
$\Phi_i$. 
Dealing with charts of polynomials implies that a coordinate system is fixed.
The signed real rational graphs will be
constructed 
 in this  coordinate system.

\subsection{$\Phi_i$ is of height 2}

Here we  consider signed real rational graphs of curves of
bidegree $(2,n)$ in $\Sigma_n$. As mentioned in section \ref{rat
  bigon}, it cannot be extracted only from the knowledge of the real
rational graphs of those curves and of their topology. However, when
such a  curve is used in a patchwork of a trigonal curve, one can recover
the sign of $a_1(X)$ out of this patchwork.

Indeed, if $\delta$ is a polygon of height $2$ in a trigonal
patchwork, let $p=(x,2)$ be the leftest point of height $2$ of $\delta$
(i.e. $x=\min_{(z,2)\in\delta}(z)$). Then, $[(0,3);p]$ is an edge of
the subdivision, and any curve of degree $3$ in the patchwork whose
Newton polygon contains this edge 
 prescribes 
 the sign of $a_1(X)$ for small
values of $X$: if the chart of one of these cubic intersect the edge
$[(0,3);p]$ (resp. $[(0,3);(-x,2)]$), then $a_1(X)$ is negative for small
positive (resp. negative) values of $X$. This is because
$[(0,3);(-x,2)]$ is of integer length $1$ and the coefficient of $Y^3$
is $1$.

For this reason, we will not write explicitly the pairs of signs on
signed real rational graphs associated to a curve of bidegree $(2,n)$
in a trigonal patchwork.

\subsubsection{$\Phi_i=\{ \delta \}$}\label{height 2}
That means that $\delta$ is the Newton polygon of a curve
$C(X,Y)=a_1(X)Y^2+a_2(X)Y+a_3(X)$ in the patchwork with $a_1(X)$
a monomial and $a_3(X)$ a non zero polynomial. Consider this curve as
a curve of bidegree $(2,n)$ in $\Sigma_n$, i.e. homogenize the curve 
$C$ in the following way :
$$C(X,Y,Z)=a_1(X,Z)Y^2+a_2(X,Z)Y+a_3(X,Z)$$
where $a_j(X,Z)$ is a homogeneous polynomial of degree $jn$. This
curve may have singular points on the fibers $\{X=0\}$ and $\{Z=0\}$.

\begin{figure}[h]
      \centering
 \begin{tabular}{ccccc}
\includegraphics[height=1.8cm, angle=0]{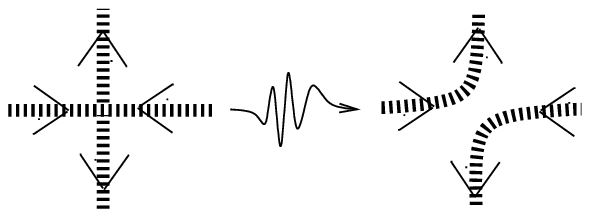}&
\hspace{5ex}&
\includegraphics[height=1.8cm, angle=0]{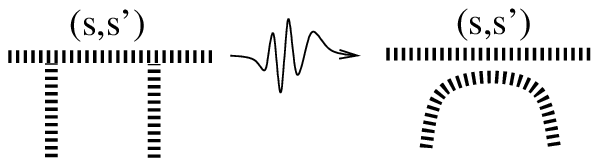}&
\hspace{5ex}&
\includegraphics[height=1.8cm, angle=0]{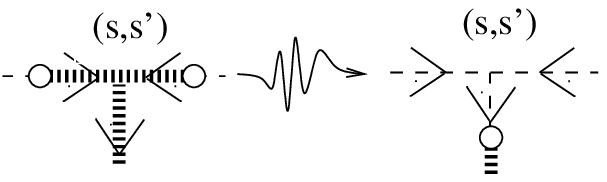}

\\ a)& &b)& &c)
\end{tabular}
\caption{}
 \label{ram pt deg 2}
\end{figure}

Consider the signed real rational graph $\Gamma_\pm(C(X,Y,Z))$
constructed in section \ref{rat bigon}. Then 
perturb this graph in such a way that 
 all the preimages of
$\infty$ 
are
 simple and perform all the operations
depicted in Figure \ref{ram pt deg 2} on $\Gamma_\pm(C(X,Y,Z))$ until
the signed real rational graph 
has no ramification point on the preimage of $[\infty,0[$. The
first operation occurs in $\mathbb CP^1\setminus\mathbb RP^1$ and the
two last ones on $\mathbb RP^1$.

The obtained signed real rational graph is  $\Gamma_\pm(\Phi_i)$.

\vspace{2ex}
\textbf{Example : }Let us detail how to find  $\Gamma_\pm(\{\delta_2 \})$ in the patchwork depicted in Figure
\ref{non convex patch}. The corresponding curve is
$C_2(X,Y)=aXY^2+bY+c$ where $a$, $b$ and $c$ are real
numbers. Considering $C_2$ as a curve in $\Sigma_4$, we have 
$$ C_2(X,Y,Z)=aXZ^3Y^2+bZ^8Y+cZ^{12}.$$
Then  $\Gamma(C_2)$ is given by the
real rational function
$$f_{\delta_2}(x)=\frac{-a^4X^4Z^{12}}{b^2Z^{16}-4acXZ^{15}-a^4X^4Z^{12}}=\frac{-a^4X^4}{b^2Z^4-4acZ^3-a^4X^4}.$$
Then, according to the chart of $C_2$, the real rational graph of
$\{\delta_2 \}$ is as depicted in Figure \ref{real graph}b). The chart
of the curve corresponding to the polygon $\delta_1$ implies that
$a_1(X)$ is negative (resp. positive) for positive (resp. negative)
values of $X$. So $\Gamma_\pm(\{\delta_2 \})$ is as depicted in
Figure \ref{real graph}b).

\subsubsection{$\Phi_i=\{ \delta_1,  \delta_2\}$
with $\delta_1$  of height 2 and  
$\delta_2$  of height 1}\label{glue 2+1}
\begin{figure}[h]
      \centering
 \begin{tabular}{ccc}
\includegraphics[height=4.3cm, angle=0]{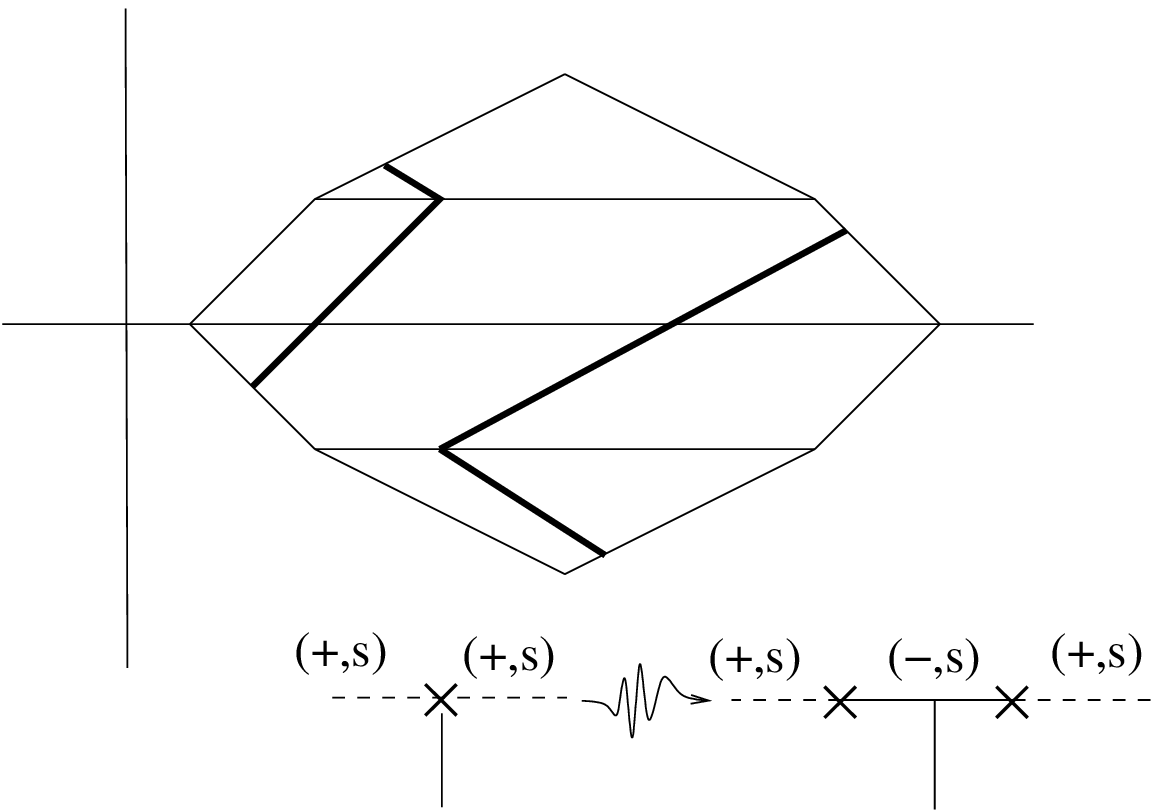}&&

\includegraphics[height=4.3cm, angle=0]{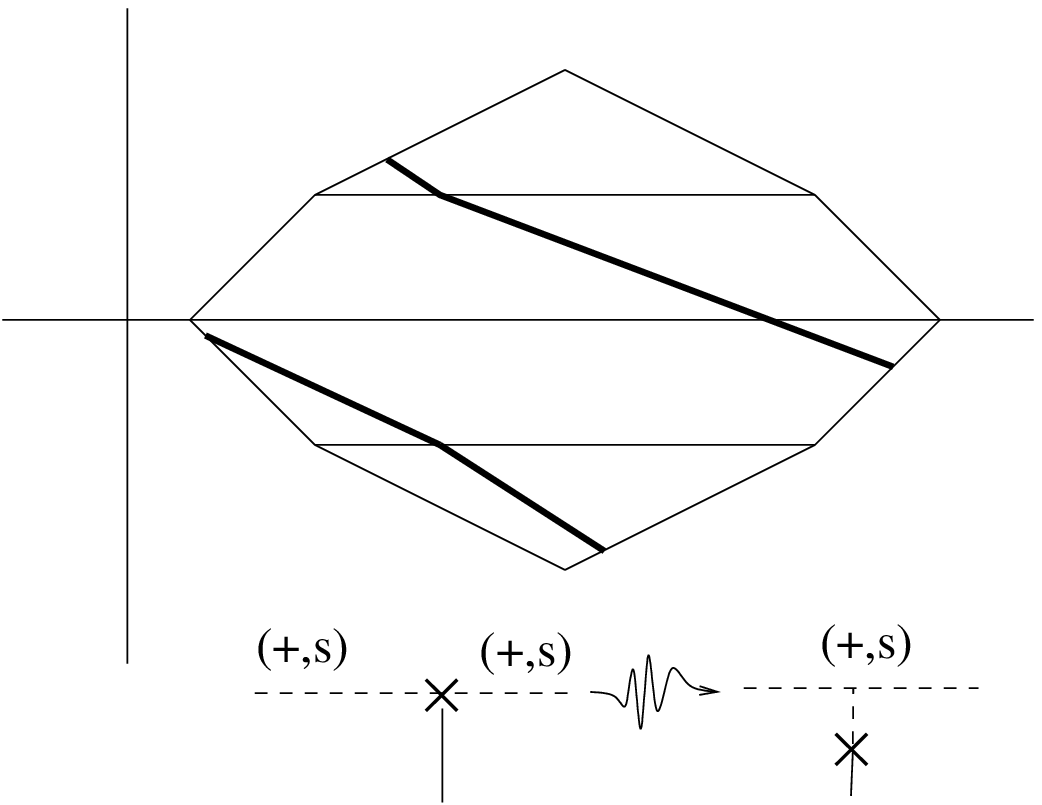}

\\ a) &&b)
\end{tabular}
\caption{}
 \label{perturb double delta}
\end{figure}

Consider the
signed real rational graph associated to the curve $C_1$ corresponding to
$\delta_1$ as described in 
section \ref{height 2}.  This curve is the
union of a curve $C_{1,1}$ of degree 1 and the curve $\{Y=0\}$. 
So
all the roots of its
 discriminant are double and
 correspond to the intersection points of $ C_{1,1}$ and
 $\{Y=0\}$. The curve $C_2$
corresponding to $\delta_2$ 
 tells us how to smooth these
 double points. More precisely, write
 $C_1(X,Y)=a_1(X)Y^2+a_2(X)Y$ and
 $C_2(X,Y)=a_2(X)Y+a_3(X)$ and suppose that $x_0$ is a
 root of $a_2(X)$ in $\mathbb R^*$. Then if $a_1(x_0)$ and $a_3(x_0)$
 have the same sign (resp. opposite signs), then perturb the
 corresponding double root of the discriminant on  $\Gamma_\pm(C_{1})$ into 2
 distinct real roots (resp. non-real roots) as depicted in Figure
 \ref{perturb double delta}a) (resp. \ref{perturb double delta}b)).

The graph  $\Gamma_\pm(\Phi_i)$  is the obtained signed rational graph.

\vspace{2ex}
\textbf{Remark : }The graph $\Gamma_\pm(\Phi_i)$ is
$\Gamma_\pm(\widetilde C)$, where $\widetilde C$ is the perturbation
of $C_{1,1}\cup\{Y=0\}$ prescribed by $ C_2$.

\subsection{$\Phi_i$ is of height 3}
\subsubsection{$\Phi_i=\{ \delta\}$}\label{height 3}

The curve $C$ in the patchwork corresponding to
$\delta$  is a real algebraic trigonal curve, nonsingular
in $\mathbb C^*\times \mathbb C$ and 
we define $\Gamma_\pm(\Phi_i)$ 
to be $\Gamma_\pm(C)$.

\vspace{2ex}
\textbf{Example : }In Figure \ref{real graph}a) (resp.\ref{real
  graph 2}b) and \ref{real graph 2}c)),
   we have depicted $\Gamma_\pm(\{\delta_1
\})$ (resp. $\Gamma_\pm(\{\delta_9
\})$  and $\{\Gamma_\pm(\delta_{10}
\})$) in the patchwork depicted in Figure
\ref{non convex patch}.

\subsubsection{$\Phi_i=\{ \delta_1,  \delta_2\}$
with $\delta_1$  of height 3 and  
$\delta_2$  of height 1}

Consider the
signed real rational graph associated to the curve $ C_1$ corresponding to $\delta_1$ as described in
section \ref{height 3}.  This curve is the
union of a curve $ C_{1,1}$ (nonsingular in $\mathbb R^*\times\mathbb R$)
of degree 2 and the curve $\{Y=0\}$ and
all the double roots of its
 discriminant 
 correspond the intersection points of $ C_{1,1}$ and
 $\{Y=0\}$. 
The curve $ C_2$ corresponding to $\delta_2$ tells us how to smooth these
 double points. More precisely, write
 $C_1(X,Y)=Y^3+a_1(X)Y^2+a_2(X)Y$ and
 $C_2(X,Y)=a_2(X)Y+a_3(X)$ and suppose that $x_0$ is a
 root of $a_2(X)$ in $\mathbb R^*$. Then if $a_1(x_0)$ and $a_3(x_0)$
 have the same signs (resp. opposite signs), then perturb the
 corresponding double root of the discriminant of the curve into 2
 distinct real roots (resp. non-real roots) as in section \ref{glue 2+1}.

The graph  $\Gamma_\pm(\Phi_i)$  is the obtained signed rational graph.

\vspace{2ex}
\textbf{Remark : }The graph $\Gamma_\pm(\Phi_i)$ is
$\Gamma_\pm(\widetilde C)$, where $\widetilde C$ is the perturbation
of $C_{1,1}\cup\{Y=0\}$ prescribed by $ C_2$. In particular,
$\Gamma_\pm(\Phi_i)$ is always the signed real rational graph associated
to a trigonal curve.

\subsubsection{$\Phi_i=\{ \delta_1,  \delta_2\}$
with $\delta_1$  of height 3 and  
$\delta_2$  of height 2}\label{glue 3+2}

 \begin{figure}[h]
      \centering
 \begin{tabular}{cc}
\includegraphics[height=3.2cm, angle=0]{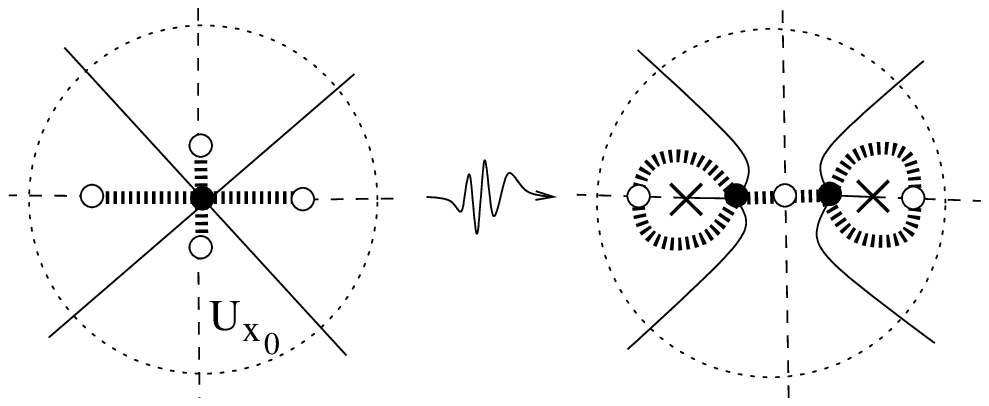}&

\includegraphics[height=2.1cm, angle=0]{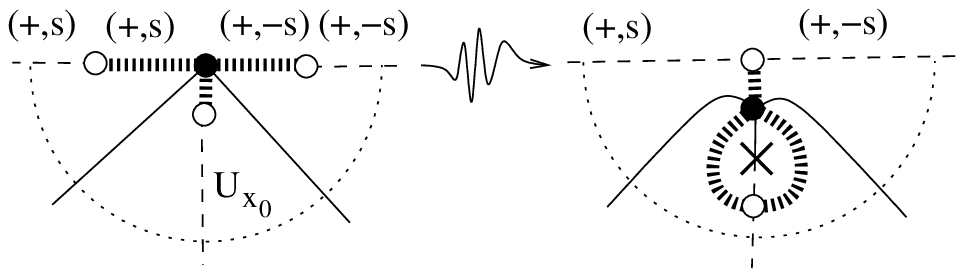}

\\ a) &b)
\end{tabular}
\caption{}
 \label{perturb a1}
\end{figure}
 \begin{figure}[h]
      \centering
 \begin{tabular}{c}
\includegraphics[height=3cm, angle=0]{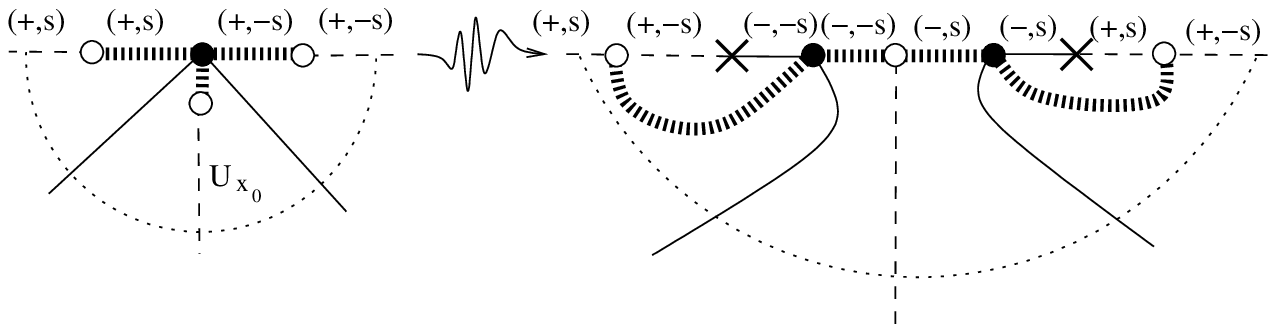}

\end{tabular}
\caption{}
 \label{perturb a1 2}
\end{figure}

First, consider  $\Gamma_\pm( \delta_2)$ as explained in section
  \ref{height 2}. We will obtain $\Gamma_\pm( \Phi_i)$
  performing some modifications on $\Gamma_\pm( \delta_2)$. 

\vspace{1ex}
Let $x_0$ be a root of $a_1(X)$ in $\mathbb C^*$.
Choose a real rational function $f$  realizing $\Gamma_\pm( \delta_2)$
 and
$\epsilon>0$ such that $f$ has no ramification point on
the connected component $U_{x_0}$ of $f^{-1}([\infty;\epsilon])$
containing $x_0$.

Then replace $U_{x_0}$ by 
\begin{itemize}
\item the part depicted in Figure \ref{perturb a1}a) if $x_0$ is non real,
\item the part depicted in Figure \ref{perturb a1}b) if $x_0$ is real
  and $a_2(x_0)<0$,
\item the part depicted in Figure \ref{perturb a1 2} if $x_0$ is real
  and  $a_2(x_0)>0$.
\end{itemize}

The graph  $\Gamma_\pm(\Phi_i)$  is the obtained signed rational graph.

\vspace{2ex}
\textbf{Example : }In Figure \ref{real graph}c), we have depicted 
 $\Gamma_\pm(\{\delta_3,\delta_4
\})$ in the patchwork depicted in Figure
\ref{non convex patch}.

\subsubsection{$\Phi_i=\{\delta_1,
 \delta_2,\delta_3\}$ 
with $\delta_1$  of height 3,  
$\delta_2$  of height 2 and $\delta_3$ of height 1}

First, consider $\Gamma_\pm(\delta_2)$ as explained in section~\ref{height 2}. Then,
perturb all the roots of $a_1(X)$ in $\mathbb C^*$ as explained in
section~\ref{glue 3+2} and perturb all the double roots of the
discriminant as explained in section~\ref{glue 2+1}.

The graph  $\Gamma_\pm(\Phi_i)$  is the obtained signed rational
graph.

\vspace{2ex}
\textbf{Example : }In Figure \ref{real graph 2}a), we have depicted 
 $\Gamma_\pm(\{\delta_6,\delta_7,\delta_8 
\})$ in the patchwork depicted in Figure
\ref{non convex patch}.

\section{Gluing of signed real rational graphs}\label{gluing graphs}

Let $\Pi$ be a trigonal patchwork of degree $n$ satisfying conditions
(1) and (2) of section \ref{state viro} and let
$\theta=(\Phi_1,\ldots,\Phi_k)$ be the $k$-tuple constructed in
section \ref{pencil}. In this section, we describe how to glue all the
signed real rational graphs $\Gamma_\pm(\Phi_i)$
 in order to obtain a new
signed real rational graph $\Gamma_\pm(\Pi)$
corresponding to a nonsingular real trigonal algebraic curve in
$\Sigma_n$. 

For $i\in\{1,\ldots k-1\}$, we describe how to glue
the signed real rational graph corresponding to $\Phi_i$ to the one
corresponding to $\Phi_{i+1}$.

\vspace{2ex}
\textbf{Example : }To illustrate this procedure, we have depicted in Figure \ref{gluing} all the
steps in the gluing process corresponding to the patchwork depicted in Figure
\ref{non convex patch}.

\subsection{$\Phi_i$ and $\Phi_{i+1}$ are of height 3 and contain no
  polygon of height 2 }

That means that both $\Phi_i$ and $\Phi_{i+1}$ contain a polygon which
have the edge $\gamma=[(0,3);(\alpha,\beta)]$ with $\beta\le1$ in
common.

As 
 noticed in section \ref{constr graph}, $\Gamma_\pm(\Phi_i)$ and
 $\Gamma_\pm(\Phi_{i+1})$ correspond to two real algebraic trigonal
 curves, respectively $C_i$ and $C_{i+1}$. In particular, the two real
 rational graphs are given by two real rational maps
$$f_i=\frac{-4P_{i}(X,Z)^3}{27Q_{i}(X,Z)^2}=\frac{R_{i}(X,Z)}{S_{i}(X,Z)}
 \mbox{ and }
 f_{i+1}=\frac{-4P_{i+1}(X,Z)^3}{27Q_{i+1}(X,Z)^2}=\frac{R_{i+1}(X,Z)}{S_{i+1}(X,Z)}.$$
 The second form of $f_i$ and $f_{i+1}$ are 
 required to be irreducible and we point out
 that
 $R_{i}(X,Z)$ can be different from $-4P_{i}(X,Z)^3$.

\vspace{1ex}
The method here is the following : identify $\infty$ in $\Gamma_\pm(\Phi_i)$ and
$0$ in $\Gamma_\pm(\Phi_{i+1})$. We obtain something which is not any more a
signed real rational graph, but 
 some
 kind of signed singular real
rational graph. This 
is a graph on two copies of $\mathbb CP^1$ which 
are identified
in one 
point, $\infty$ for one of them and $0$ for the other one. Smoothing
equivariantly this
point  we 
 will
 obtain a smooth signed
real rational graph on  $\mathbb CP^1$. 

\subsubsection{$\gamma=[(0,3);(3l,0)]$ or $\gamma=[(0,3);(2l,1)]$}\label{non sing fiber}
 From condition $(1)$ we get that
$f_i(\infty)\ne 1$ and
$f_{i+1}(0)\ne 1$. Moreover,
one can perturb slightly the coefficients corresponding to the integer
vertices of $\gamma$ such that $f_i(\infty)\in\mathbb
RP^1\setminus\{0,1,\infty\}$. This implies that $f_{i+1}(0)\in\mathbb
RP^1\setminus\{0,1,\infty\}$.

 \begin{figure}[h]
      \centering
 \begin{tabular}{ccc}
\includegraphics[height=2.5cm, angle=0]{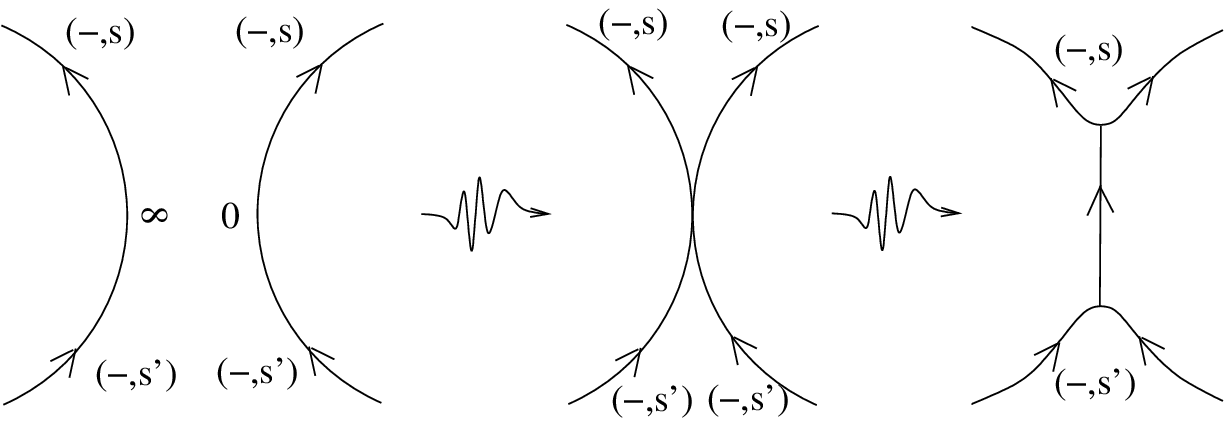}&
&\hspace{4ex}
\includegraphics[height=2.5cm, angle=0]{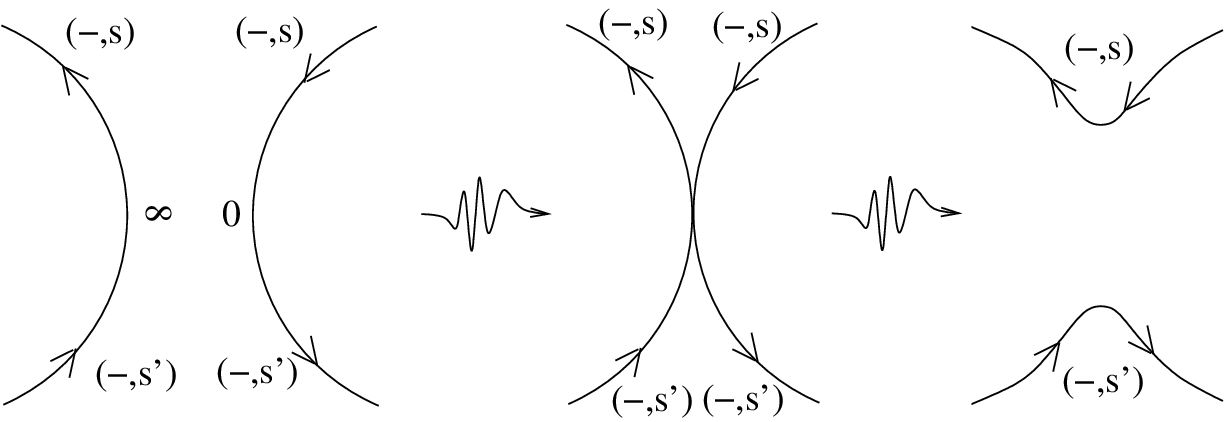}

\\ a)& &b)
\end{tabular}
\caption{}
 \label{smooth rat 1}
\end{figure}

Then, as the truncation
  on $\gamma$ of $C_i$ and $C_{i+1}$ coincide, the pair of signs on
  the positive (resp. negative) part of $\Gamma_\pm(\Phi_i)$  near $\infty$
   and $\Gamma_\pm(\Phi_{i+1})$ near 
$0$
  also coincide. Suppose that
  $f_i(\infty)$ is
  in $]0,1[$. Then, smooth the singular real rational graph as depicted in
  Figure \ref{smooth rat 1}a) or b), depending on the situation. The
  case when $f_i(\infty)$ is in $]1,\infty[$ or $]\infty,0[$ can be
  treated analogously.

 \begin{figure}[h]
      \centering
 \begin{tabular}{cc}
\includegraphics[height=3cm, angle=0]{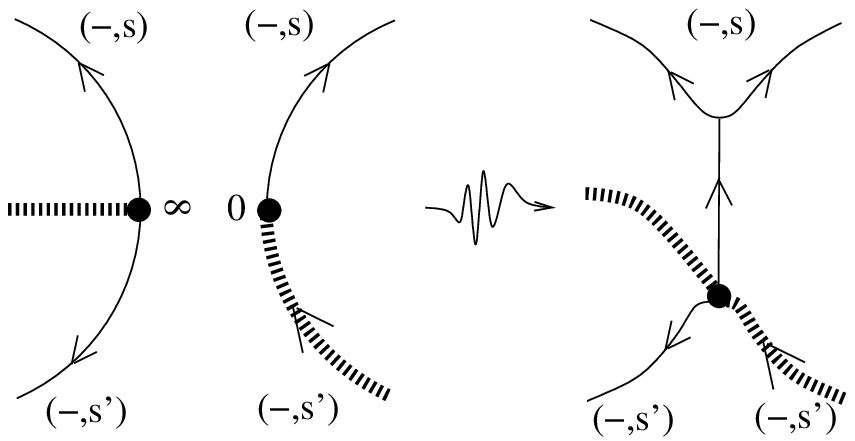}&

\\a)

\end{tabular}
\caption{}
 \label{smooth rat 6}
\end{figure}

\subsubsection{$\gamma=[(0,3),(3l-1,0)]$ or $\gamma=[(0,3),(3l-2,0)]$}
We treat the case $\gamma=[(0,3),(3l-1,0)]$, the case
$\gamma=[(0,3),(3l-2,0)]$ being symmetric.

In this case, a simple computation shows that $\infty$ is a root of
order $3a+1$ of $R_i$ and that $0$ is a root of order $3b+2$ of
$R_{i+1}$, where $a$ and $b$ are some natural numbers. Perturbing if
necessary the coefficient of $C_i$ and $C_{i+1}$, one can suppose $a=b=0$.
Moreover, as the truncation
  on $\gamma$ of $C_i$ and $C_{i+1}$ coincide, the pair of signs on
  the positive (resp. negative) part of $\Gamma_\pm(\Phi_i)$  near $\infty$
   and $\Gamma_\pm(\Phi_{i+1})$ near $0$
  also coincide.
Then, smooth the singular real rational graph as depicted in
  Figure \ref{smooth rat 6} or
symmetrically with respect to the  axis $\{Y=0\}$,
  depending on the situation.

\subsubsection{$\gamma=[(0,3),(2l-1,1)]$}

 \begin{figure}[h]
      \centering
 \begin{tabular}{ccc}
\includegraphics[height=3cm, angle=0]{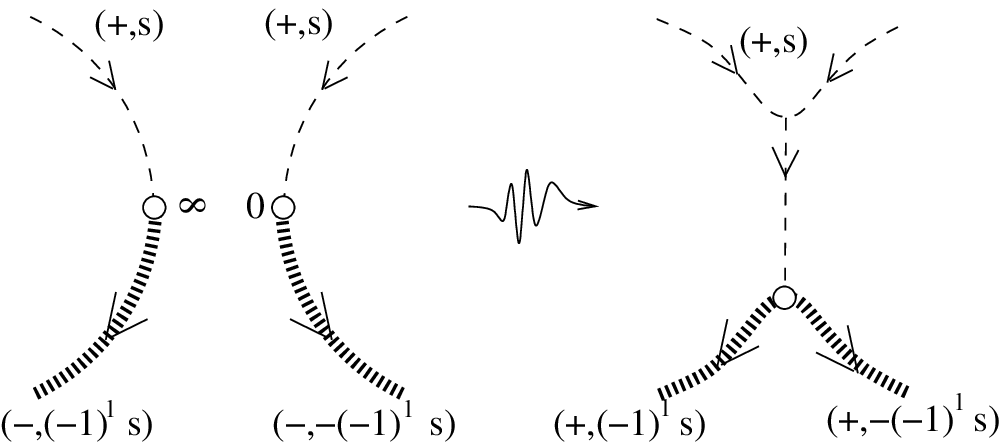}&
\hspace{4ex}&
\includegraphics[height=3cm, angle=0]{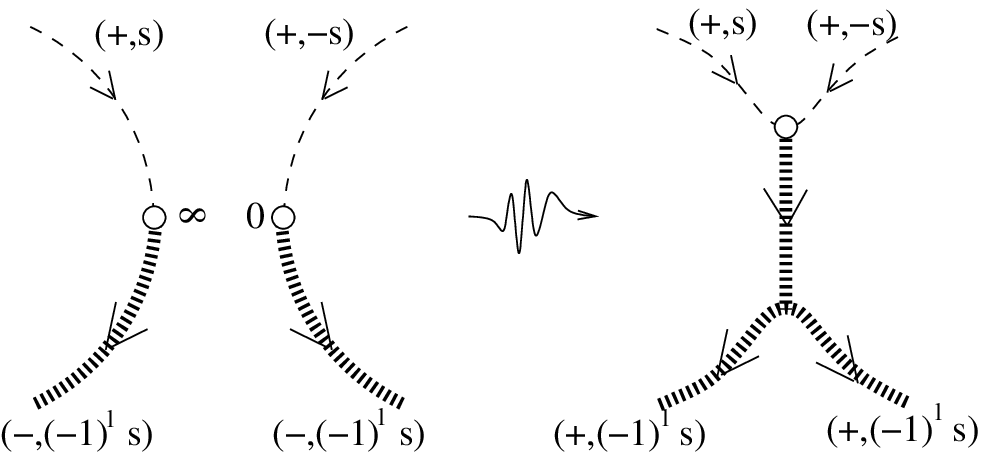}

\\a)&&b)

\end{tabular}
\caption{}
 \label{smooth rat 7}
\end{figure}

In this case, a simple computation shows that $\infty$ is a root of
order $2a+1$ of $S_i$ and that $0$ is a root of order $2b+1$ of
$S_{i+1}$, where $a$ and $b$ are some natural numbers. Perturbing if
necessary the coefficient of $C_i$ and $C_{i+1}$, one can 
assume that 
 $a=b=0$.
If
  $(s_D,s_Q)$ is the pair of signs on  
 the positive part of $\Gamma_\pm(\Phi_i)$
 (resp. $\Gamma_\pm(\Phi_{i+1})$) near $\infty$ (resp. $0$), then the pair
  of signs on  
 the negative part of $\Gamma_\pm(\Phi_i)$ (resp. $\Gamma_\pm(\Phi_{i+1})$) near $\infty$ (resp. $0$) is
  $(s_D,(-1)^{l}s_Q)$ (resp. $(s_D,-(-1)^{l}s_Q)$). 
Moreover, as the truncation
  on $\gamma$ of $C_i$ and $C_{i+1}$ coincide, the sign of the
  discriminant on
  the positive (resp. negative) part of $\Gamma_\pm(\Phi_i)$  near $\infty$
   and $\Gamma_\pm(\Phi_{i+1})$ near $0$
  also coincide.

Then, depending on the situation, smooth the singular real rational
graph using the appropriate perturbation among those depicted in
  Figure \ref{smooth rat 7} and  their symmetric with respect to the
  axis $\{Y=0\}$.

\subsection{$\Phi_i$ or $\Phi_{i+1}$  contains a
  polygon of height 2}

\begin{figure}[h]
      \centering
 \begin{tabular}{cc}
\includegraphics[height=5cm, angle=0]{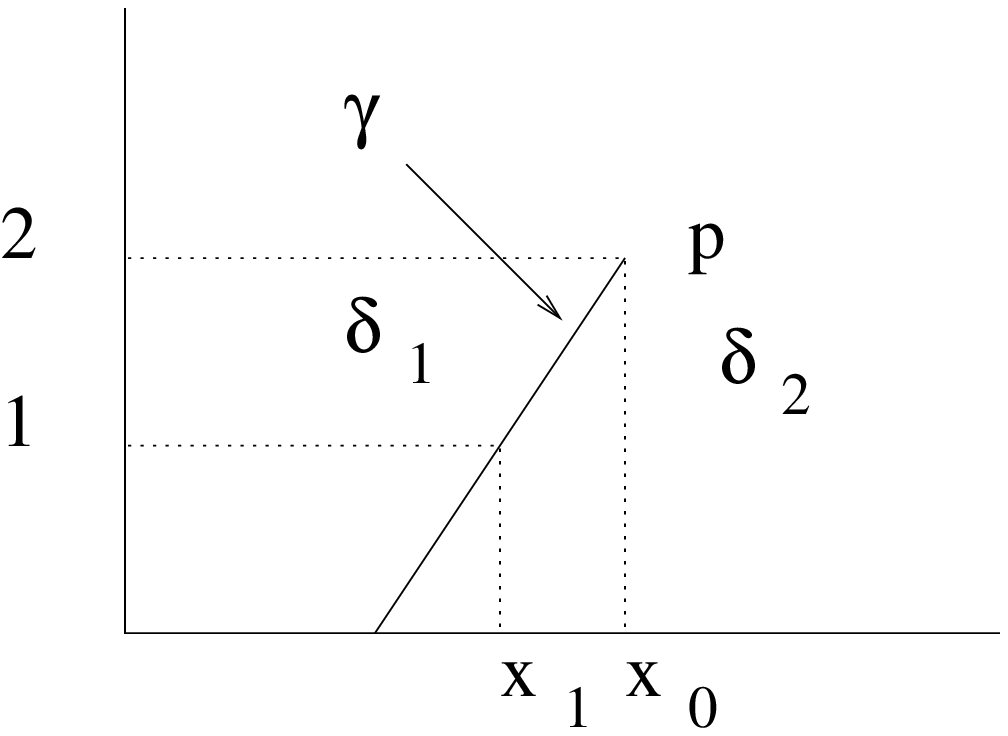}&
\includegraphics[height=5cm, angle=0]{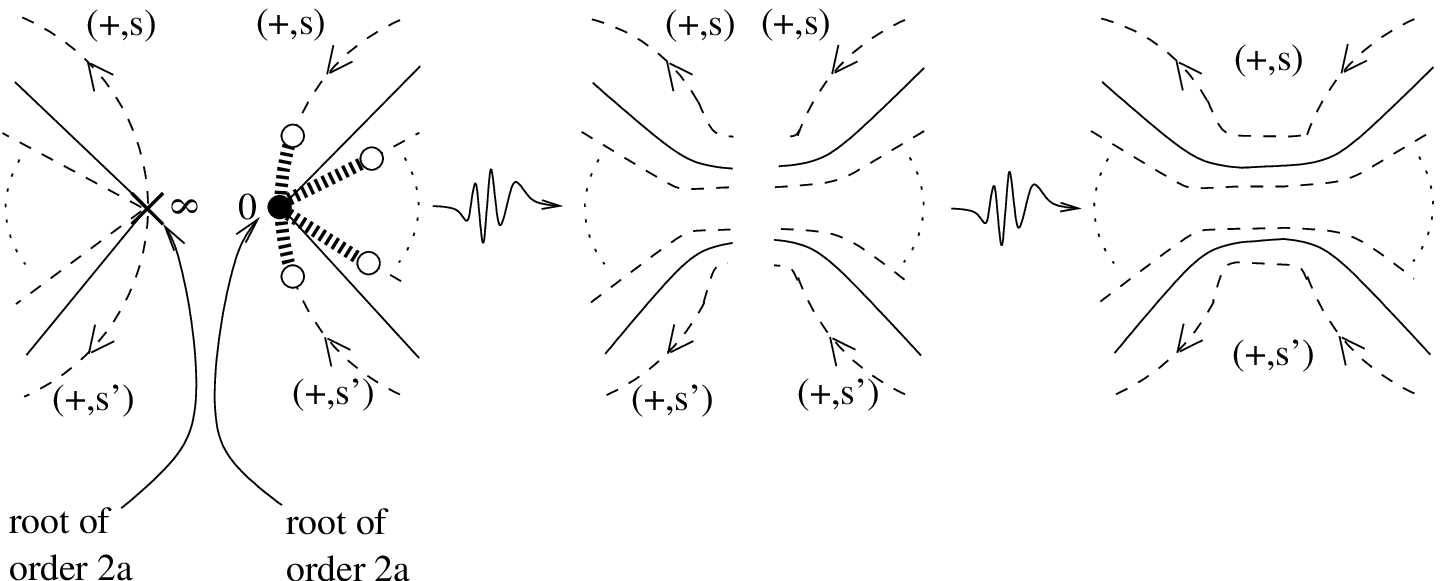}

\\a)&b)

\end{tabular}
\caption{}
 \label{smooth rat 8}
\end{figure}

Let $p=(x_0,y_0)$ be the point of the polygons of $\Phi_i$ of ordinate 2 with
maximal abscissa. Let $\gamma$ be the unique edge of the form $[p;(\alpha,\beta)]$ with $\beta<2$
which is a common edge of a polygon $\delta_i$ of
 $\Phi_i$ and of a polygon $\delta_{i+1}$ of $\Phi_{i+1}$.
The line supported by $\gamma$ intersects the line $\{Y=1\}$ in the
point  $(x_1,1)$ (see
figure \ref{smooth rat 8}).

We suppose that $x_1\le x_0$. The case  $x_1\ge x_0$ can be recovered by
symmetry.

\vspace{2ex}
Let $C_i$ (resp. $C_{i+1}$) be the curve of the patchwork
corresponding to $\delta_i$ (resp. $\delta_{i+1}$) and  $f_i=\frac{R_{i}(X,Z)}{S_{i}(X,Z)}$
(resp. $f_{i+1}=\frac{R_{i+1}(X,Z)}{S_{i+1}(X,Z)}$) be the real rational map constructed out of $C_i$
(resp. $C_{i+1}$)  as explained in section \ref{graph}. Here again,
the polynomials $R_{i}(X,Z)$ and $S_{i}(X,Z)$ (resp.$R_{i+1}(X,Z)$ and
$S_{i+1}(X,Z)$) have no common factor.

According to section \ref{constr graph}, the graph $\Gamma_\pm(\Phi_i)$
(resp. $\Gamma_\pm(\Phi_{i+1})$) is either $\Gamma_\pm(C_i)$
(resp. $\Gamma_\pm(C_{i+1})$) or a perturbation of this signed real
rational graph in some points distinct from $0$ and $\infty$.
So, the gluing of $\Gamma_\pm(\Phi_i)$
and $\Gamma_\pm(\Phi_{i+1})$ is determined by the neighborhood of
$\infty$ in $\Gamma_\pm(C_i)$ and $0$ in  $\Gamma_\pm(C_{i+1})$.

A simple calculation shows that
\begin{itemize}
\item if $x_1$ is an integer, then $\infty$ (resp. $0$) is a root of
  order $2(2x_0-x_1)$ of $R_i-S_i$ (resp. $R_{i+1}$),

\item if $x_1$ is not an integer, then $\infty$ (resp. $0$) is a root of
  order $2(2x_0-[x_1])-1$ of $R_i-S_i$ (resp. $R_{i+1}$),

\end{itemize}
where $[x]$ denote the integer part of the real $x$.

It follows from the construction of $\Gamma_\pm(\Phi_i)$ and
  $\Gamma_\pm(\Phi_{i+1})$ and from condition (1) of section \ref{state viro}
  that the pair of signs on 
  the positive (resp. negative) part of $\Gamma_\pm(\Phi_i)$  near $\infty$
   and $\Gamma_\pm(\Phi_{i+1})$ near $0$
  also coincide.

Let 
$U_0$ be the connected component of $f^{-1}_{i+1}([\infty,0])$ which contains $0$.
Then cut $U_0$ from $\Gamma_{i+1}$, cut $\infty$ from $\Gamma_i$ and
glue it in the only possible way, as it is depicted on Figure \ref{smooth
  rat 8}b) in a particular case.

 \begin{figure}[!h]
      \centering
 \begin{tabular}{ccccc}
\includegraphics[height=3.1cm, angle=90]{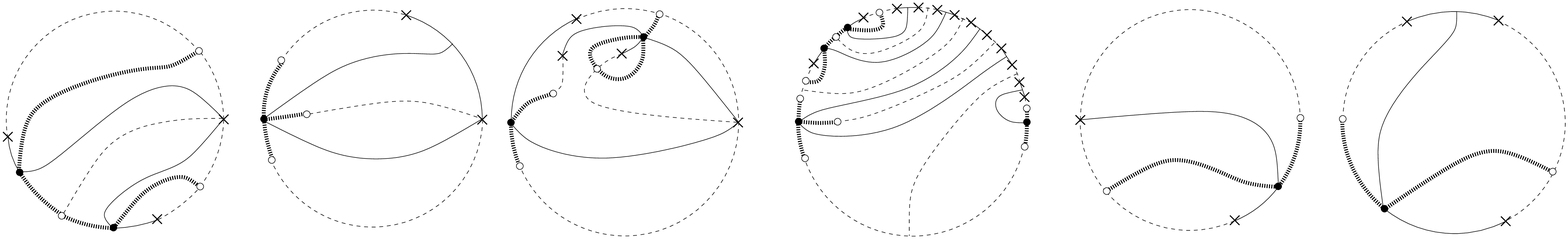}&
\hspace{10ex} &
\includegraphics[height=3.1cm, angle=90]{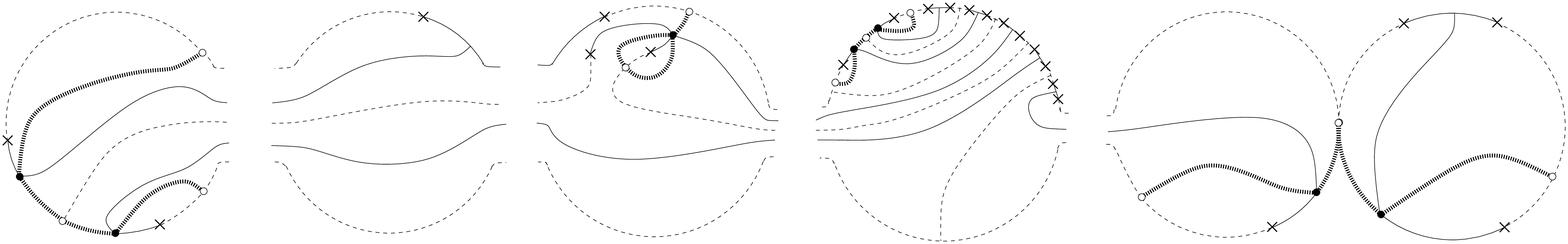}&
\hspace{10ex} &

\includegraphics[height=3.1cm, angle=90]{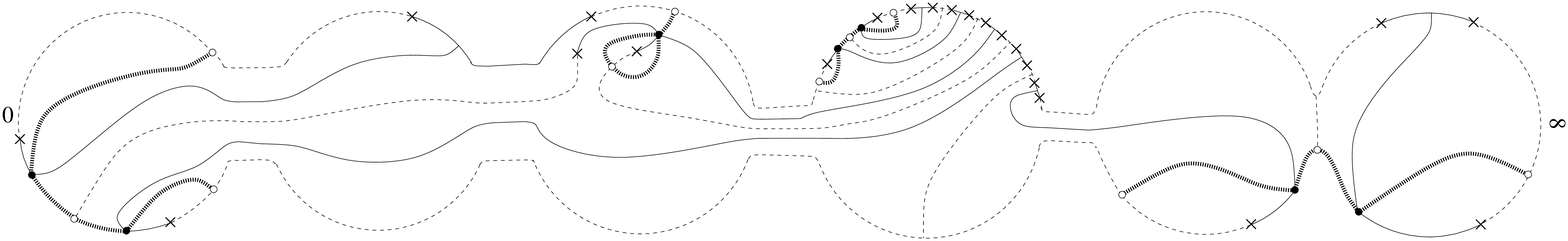}

\\ a) &&b)&&c)
\end{tabular}
\caption{}
 \label{gluing}
\end{figure}

\section{Proof of Theorem \ref{non convex viro}}\label{proof}

Let $\Pi$ be a trigonal patchwork of degree $n$ satisfying conditions
(1) and (2) of section \ref{state viro}.
In this section, we show that the signed real rational graph $\Gamma_\pm(\Pi)$
 constructed in section \ref{gluing graphs} out of $\Pi$  corresponds to a
$\mathcal L$-nonsingular real trigonal algebraic curve in
$\Sigma_n$ which realizes the sign array of the corresponding patchwork.

\vspace{2ex}
First, we check that $\Gamma_\pm(\Pi)$ is the signed real rational
graph of a $\mathcal
L$-nonsingular real trigonal algebraic curve in $\Sigma_n$.
The only non-trivial part is to prove that the degree of
 $\Gamma_\pm(\Pi)$ is  $6n$.  
 In order to compute it, we count the number of
 preimages of 1. So we first have to count how many preimages
 of 1 
 each element of $\theta$ brings to  $\Gamma_\pm(\Pi)$. According to
 section \ref{gluing graphs}, 
this is equivalent to count how many preimages of $1$ are lying on
$\Gamma(\Phi_i)\setminus\{0,\infty\}$ for all element $\Phi_i$ of
$\theta$.

First, we count the number of roots in $\mathbb C^*$ of the
discriminant of a curve.
Given a convex polygon $\delta$, the quantity $I(\delta)$ is the number
of integer points in the interior of $\delta$ and $\partial(\delta)$
denotes the number of integer points of the boundary of $\delta$ which
have neither a maximal nor a minimal ordinate.

\begin{lemma}\label{pt tg}
Let $C(X,Y)$ be a totally nondegenerate curve  with Newton polygon
$\delta$. Then, the discriminant of
$C(X,Y)$ has $2I(\delta)+\partial(\delta)$ roots in $\mathbb C^*$
counted with multiplicity.
\end{lemma}
\textit{Proof : }Let us consider the closure of $C(X,Y)$ in $\mathbb
CP^2$ and the rational map $\pi : \mathbb CP^2 \dashrightarrow \mathbb
CP^1$ given by $\pi([x:y:z])=[x:y]$. It is well known that the genus
of the normalization of $C$ is $I(\delta)$. So, applying the Riemann
Hurwitz formula to the restriction of $\pi$ to $C$, we have

$$2-2I(\delta) = 2d-r_1-r_2 $$

where $d=\deg_Y(C)$, $r_1$ is the number of roots (counted with
multiplicity) in $\mathbb C^*$ of the discriminant of $C(X,Y)$ and
$r_2$ is the sum of the tangency order of
 local 
branches
 of $C$ 
along
$\{X=0\}$ and $\{Z=0\}$. As $C(X,Y)$ is totally nondegenerate, we
have $r_2=2d-\partial(\delta)-2$, so

$$r_1=2I(\delta)+\partial(\delta) .$$

The Lemma is proved.\findemo

\begin{cor}\label{deg gam}
The real rational graph $\Gamma_\pm(\Pi)$ constructed out of a
trigonal patchwork of degree $n$ is
of degree $6n$.
\end{cor}
\textit{Proof : }It is clear from the construction and the hypothesis
on the trigonal patchwork that the degree of $\Gamma_\pm(\Pi)$ is equal to the
sum of the numbers of preimages of $1$ in $\Gamma(\Phi_i)\setminus \{0,\infty\}$ for all
$\Phi_i$ in $\theta$. 
 From
 section \ref{constr graph} and Lemma \ref{pt
  tg},
 we know that
 this number is equal to $2I(\Delta)+\partial(\Delta)$ which is 
equal to $6n$.\findemo

\begin{prop}\label{almost done}
The signed real rational graph $\Gamma_\pm(\Pi)$ is realizable by a $\mathcal
L$-nonsingular real algebraic trigonal curve $C$ in $\Sigma_n$. 
\end{prop}
\textit{Proof : }
 The degree of $\Gamma_\pm(\Pi)$  is given by Corollary \ref{deg gam} and
 all the other conditions of Proposition \ref{iff trigo} are
 satisfied by construction. \findemo

Choose a standard coordinate system on $\Sigma_n$ such that the point
$\infty$ on $\Gamma_\pm(\Pi)$ corresponds to the fiber at infinity and
let $SA_C$ be the sign array of the curve $C$ in this coordinate system.

\begin{prop}\label{non convex viro prop}
Under the hypothesis (1) and (2) of section \ref{state viro}, the sign
array $SA_C$ and the sign array constructed out of $\Pi$ in section
\ref{topology trigonal} coincide.
\end{prop}
\textit{Proof : }Let $[s_0,s_1\ldots s_u]$ (resp. $[t_0,t_1\ldots
  t_v]$) be the sign array $SA_C$ (resp. associated to the
  patchwork).
First note that a direct consequence of sections  \ref{sign
  patchwork}, \ref{constr graph} and
  \ref{gluing graphs}
  is that $u=v$, i.e. both curves have the same
  number of tangency points with the pencil of vertical lines. Hence,
  to prove Theorem \ref{non convex viro}, we have just have to prove
  that
$s_j=t_j$ for any $j$.

It is clear from sections \ref{sign patchwork}, \ref{constr graph} and \ref{gluing graphs} 
   that both $s_j$ and $t_j$ come from the same
  real $x_j$ for some $\Phi_i$.  Choose a real rational map
   $f=\frac{-4P^3}{27Q^2}$ which 
   realize $\Gamma_\pm(\Pi)$.

If $\Phi_i$ contains no polygon with a horizontal edge of 
height
 $2$, then $x_j$
corresponds to a root of the discriminant of a curve of $\Pi$ whose Newton
polygon is in $\Phi_i$. Hence,  $s_j$ which is the sign of $Q(x_j)$, is
the second element of the pair of signs which labels
$\Gamma_\pm(\Phi_i)$ near $x_j$. But according to section \ref{sign
  patchwork}, this is also $t_j$.

If $\Phi_i$ contains a polygon with a horizontal edge of 
height $2$, denote by $a_1(X)$ the polynomial corresponding to this
edge. Let $a$ be
the number of roots of $a_1(X)$  counted with multiplicity which are
strictly between $0$ and $x_j$ and let $\sigma$ be  the sign of
  $a_1(x)$ for
$x\in\mathbb R^{*}$ small enough with the same sign as $x_j$.
There are two possibilities, each of them can checked directly from
the construction of 
$\Gamma_\pm(\Phi_i)$. 
\begin{itemize}
\item If $x_j$ is not a root of $a_1(x)$, then  $s_j=(-1)^a\sigma=t_j$.
\item If $x_j$ is a root of $a_1(x)$, then $x_j$ gives rise to two
  elements of both sign arrays, say $s_j$, $s_{j+1}$, $t_j$ and
  $t_{j+1}$. We have $s_j=-(-1)^a\sigma=t_j$ and $s_{j+1}=(-1)^a\sigma=t_{j+1}$.

\end{itemize}
So in any case, we have $s_j=t_j$.~\findemo

\vspace{2ex}
\textbf{Example : }One can check  on Figure \ref{gluing}c), that the
obtained real rational graph corresponds to a non singular real
algebraic trigonal curve realizing the sign array of 
 the patchwork depicted in Figure
\ref{non convex patch}.

\section{Remark on an equivariant version Theorem \ref{viro}}

We were guided in this work by the following observation : the
position of the real part of a nonsingular real algebraic curve in
$\Sigma_n$ with respect to the pencil of lines is 
totally encoded by the sign array realized by some collection of
polynomials in one variable, namely the leading coefficients of some
signed subresultant sequences (see \cite{BPR} or \cite{Hon2} for
example). It seems to us that this observation could lead to
interesting results in real algebraic geometry, especially in the
study of the need of convexity 
in patchworking.

In the case of trigonal curves, the polynomials whose sign array
encode the topology of the real part of the curve are exactly, with notation
of section \ref{poly deg 3}, the polynomials $P(X)$, $Q(X)$ and
$D(X)$. Moreover, such three polynomials with a prescribed sign array
can be constructed using 
signed real rational graphs.

\vspace{2ex}
On the other hand, real rational graphs carry much more informations than
the knowledge of the topology of the real part of real algebraic
curves. Indeed, one can extract  
the position with respect to the pencil of lines of the
complexification of any real trigonal curve out of its associated real
rational graph. Hence, one can consider not only isotopies
of curves in the real part of the surface $\Sigma_n$ but also
equivariant isotopies of curves in $\Sigma_n$. That means that it is
possible to state an equivariant version of Theorem \ref{viro} and
Theorem \ref{non convex viro}.

Proving such an equivariant statement would carry us quite
far away from our original framework. So we only give  the main ingredients
needed in its formulation and in its 
proof.

\begin{description}
\item[Step 1. ]In the patchwork construction, one can use
  complexification of lattice polygons and 
complex
  equivariant charts of real 
  polynomials as defined in \cite{IS}, \cite{IS2} and \cite{V2}. Then,
  the patchwork procedure gives a piecewise smooth surface,
  invariant under the action of the complex conjugation, in the
  complexification of the triangle with vertices $(0,0)$, $(0,3)$ and $(3n,0)$.

\item[Step 2. ]One can use the braid monodromy factorization (see
  \cite{KM1} or \cite{KK1} for example) to
  encode the position of a 
  complex algebraic curve in $\Sigma_n$ with respect to the pencil of
  lines. An equivariant version of this encoding, based on previous
  papers by Orevkov (see \cite{O1} for example), should be used to deal
  with real algebraic curves.

\item[Step 3. ]As we associated a sign array to a trigonal patchwork
  (in section \ref{sign
  patchwork}), one can associate in a similar way  an
  equivariant braid monodromy 
  factorization to an equivariant trigonal patchwork (defined in Step
  1).

\item[Step 4. ] One can extract the equivariant braid monodromy
  factorization 
  of the real trigonal algebraic curve constructed in section
  \ref{proof} out of the signed real rational graph
  constructed in section  \ref{gluing graphs} and prove that it
  coincides with the one constructed in Step 3.

\end{description}

\small
\def\rightmark{\em Bibliography}
\addcontentsline{toc}{section}{References}
\bibliographystyle{alpha}
\bibliography{/home/erwan/Maths/Recherche/Biblio}

\vspace*{2 ex}
\begin{tabular}{lll}
\textbf{Beno\^{i}t Bertrand}& \hspace{10ex}  & \textbf{Erwan Brugallé}\\
Section de math\'ematiques &    &          Max Planck Institute für Mathematik\\
Universit\'e de Gen\`eve &    &  Vivatsgasse, 7\\
 case postale 64, &          &                 D-53111 Bonn \\
2-4 rue du Li\`evre,          &             &                  Deutschland\\
Gen\`eve &&\\
Suisse &&\\
\\
E-mail : benoit.bertrand@math.unige.ch & &E-mail : brugalle@mpim-bonn.mpg.de

\end{tabular}

\end{document}